\documentclass[pdftex,a4paper]{beamer}
\mode<presentation>
{
\usetheme[secheader]{Boadilla}
}
\usepackage[english]{babel}
\usepackage{mathrsfs}
\usepackage{amsmath}
\usepackage{ams}
\usepackage{color}
\newcommand{\blu}[1]{{\textcolor{blue}{#1}}}
\newcommand{\pur}[1]{{\textcolor{purple}{#1}}}
\newcommand{\nz}{\normalsize}\newcommand{\sz}{\small}
\newcommand{\fz}{\footnotesize}
\newcommand{\Ds}{\displaystyle} \newcommand{\Ts}{\textstyle}
  
\newlength{\ei}\ei=0.0138888889em   
\newcommand{\gP}{{\mbox{\Large.}}}
\newcommand{\so}{\;\Longrightarrow\;}
\newcommand{\Lo}{\mathop{\sf {{}o{}}}\nolimits}
\newcommand{\Ew}{\mathop{\rm {}E{}}\nolimits}
\newcommand{\Cov}{\mathop{\rm {{}Cov{}}}\nolimits}
\newcommand{\lin}{\mathop{\rm lin}\nolimits}
\newcommand{\cl}{\mathop{\it c\ell}\nolimits}
\newcommand{\clin}{\mathop{\cl\lin}\nolimits}
\newcommand{\MSE}{\mathop{\rm MSE}\nolimits}
\newcommand{\EM}{{\mathbb I}}
\newcommand{\Jc}{\mathop{\sf 1}\nolimits}
\newcommand{\lJc}{\mathop{\large\sf 1}\nolimits}
\newcommand{\Lw}{\mathop{\mathsurround0em\mbox{${\cal L}$}}\nolimits}
\newcommand{\diag}{\mathop{\sf diag}\nolimits}
\newcommand{\Nlim}{\lim\nolimits}

\newcommand{\Nsup}{\sup\nolimits}

\newcommand{\Tfrac}[2]{{\Ts\frac{#1}{#2}}}
\newcommand{\Tsum}{\mathop{\Ts\sum}\nolimits}
\newcommand{\Tint}{\mathop{\Ts\int}\nolimits}
\newcommand{\maxev}{\mathop{\rm max\hspace*{6\ei}ev}\nolimits}

\newcommand{\zi}{\mathop{\sf sign}\nolimits}
\newcommand{\sign}{\mathop{\sf sign}\nolimits}
\newcommand{\tr}{\mathop{\sf tr}\nolimits}

\newcommand{\med}{\mathop{\sf med}\nolimits}
\title[Semiparametrics and Robustness]{Connections between Semiparametrics and Robustness}
\author[Helmut Rieder]{Helmut Rieder}
\institute[UBT]{University of Bayreuth, Germany}
\date[23/05/2013]{Tsinghua University, Beijing, 23 and 28 May 2013}
\begin{document}
\begin{frame}  \titlepage \end{frame}
\section{Introduction}
\begin{frame}\null
\blu{\nz Topics} \nz
\par \bigskip
1) Robust influence curves for models 
with $\infty$-dim.\ nuisance parameter;
e.g., semiparametric regression (Cox), mixture models (Neyman--Scott).
\par\smallskip
2) Adaptiveness (Stein's necessary condition) of robust 
estimators w.r.t.\ a finite-dim.\ nuisance parameter;
e.g., location, linear regression, and ARMA.
\par \smallskip
3) Semiparametric treatment of gross error deviations from an ideal model as an
   $\infty$-dim.\ nuisance parameter, by projection on balls; for testing,
   an asymptotic version of the Huber--Strassen maximin result is thus obtained.
\par \smallskip
4) Uniform and nonuniform asymptotic normality of robust and adaptive estimators, respectively,
in regression and time series models.
\smallskip \par
5) Fragility of optimal one-sided tests and confidence limits obtained for convex tangent cones, by
projection on cones, as opposed to stability of corresponding procedures, even two-sided,
for linear tangent spaces.
\smallskip \par
6) Control of the unknown neighborhood radius, a nuisance parameter in robustness.
\par \null \vfill
\end{frame}
\section{1. Semiparametric Setup}
\subsection{1.1 Model, ICs, functionals, ALEs}
\begin{frame} 
\blu{Model} $ {\cal Q}= \{\,Q _{\theta,\nu}\,\}$ with parameter of interest
   $\theta\in \Theta$ open $\subset$ some $ \R^k$, and
   nuisance parameter $\nu\in H_\theta$.
\blu{Differentiability} at any fixed $(\theta_0, \nu_0)$
\begin{equation}
   dQ_{\theta_0+ta, \nu^g_t} \approx
      \bigl(1+ t (a' \Lambda+g)\bigr)
      dQ_{\theta_0,\nu_0} \quad\mbox{as $t\to0$}
\end{equation}
in direction~$a\in\R^k$, along paths~$t \mapsto\nu^g_t$.
\blu{Tangents}\enskip $g\in L_2(Q _{\theta_0,\nu_0})$, $g \perp1$, 
   $$\partial_1 {\cal Q}=\{\,a' \Lambda\mid a\in\R^k\,\}\,,\quad
        \partial_2 {\cal Q}\hskip36\ei\mbox{a cone,}\quad
   \partial {\cal Q}= \partial_1 {\cal Q} + \partial_2 {\cal Q}
   $$\\
\blu{Fisher information} of ${\cal Q} _{\nu_0}$ ($\nu$ fixed to~$\nu_0$)
about $\theta$ at $\theta_0$: \quad ${\cal I}=\Cov \Lambda$
\par \hfill\mbox{\tiny Bickel (1982), Bickel et al. (1993), v.d.Vaart (1998)}
\par \vspace{-\medskipamount}%
\blu{Influence curves} at $Q _{\theta_0,\nu_0}$:
                         \hfill\mbox{\tiny Rieder (1994), Shen (1995)}
\vspace{-.5 \smallskipamount}\begin{equation}
  \hspace{-2em}
  \psi \in L_2^k,\quad \Ew \psi=0, \quad \Ew \psi \Lambda'= \EM_k,
  \quad \Ew \psi g =0\enskip \forall g\in \partial_2 {\cal Q}
\vspace{-.5 \smallskipamount}\end{equation}
\blu{F-consistent diff.\ functionals:}\quad
   $    T(Q_{\theta_0+ta,\nu^g_t}) - \theta_0 \approx
        \Ew \psi (a' \Lambda+g)\,t         = t a $
\par \smallskip
\blu{AL estimators:}\hfill
  $    n ^{1/2}(S_n- \theta_0) \approx n ^{-1/2}\sum_{i=1}^n \psi(x_i) $\quad
in $ Q _{\theta_0,\nu_0} ^{(n)}$-probability,\\
such that\vspace{-\medskipamount}
   $$    \sqrt{n}\,(S_n- \theta_0)\longrightarrow {\cal N}(a, \Cov\psi)
         \quad \mbox{in law under~$Q_n ^{(n)}(a,g)$} \hspace{-5em}
   $$
where $Q_n(a,g)=Q_{\theta_0+ s_na,\,\nu_{s_n}^{g}} $, at scale $s_n=1/\!\sqrt{n}\,$,
      $n=$ sample size.
\end{frame}
\subsection{1.2 Efficient IC, classical adaptivity}
\begin{frame}
Let $\Pi$, $\Pi_2 $ denote the (coordinatewise) orthogonal projections
from $L_2^k(Q _{\theta_0,\nu_0}) $ on the closed linear spans
  $ \clin \partial {\cal Q} = \partial_1 {\cal Q} + \clin \partial_2 {\cal Q}$
and $\clin \partial_2 {\cal Q}$, respectively.
\blu{Unique projection on $\clin \partial {\cal Q}$} of all ICs: 
\begin{equation}
   \Pi(\psi) = \blu{ \psi _{\sf eff}:
   = {\cal J}^{-1} \bar{\Lambda} }\qquad \forall\,\psi \: \mbox{IC}
\end{equation}
where $ \bar{\Lambda}:=\Lambda- \Pi_2(\Lambda)$ (model ${\cal Q}$).
  \blu{ $\psi _{\sf class}:={\cal I}^{-1}\Lambda $ }(model~${\cal Q} _{\nu_0}$). \\
Fisher informations of~${\cal Q}$ and~${\cal Q} _{\nu_0}$ for~$\theta$
at~$(\theta_0,\nu_0)$ and~$\theta_0$, respectively:
\begin{equation}
   {\cal J} = \Cov\bar{\Lambda}  =
               {\cal I} - \Cov \Pi_2(\Lambda) \le {\cal I}=\Cov \Lambda
\end{equation}
\blu{Asymptotic covariance bound} for AL estimators with IC~$\psi$:
\begin{equation}
     \Cov\psi \ge  {\cal J}^{-1}=\Cov\psi _{\sf eff} \,, \hspace{1.5em}
     \mbox{attained iff}\hspace{.75em} \psi=\psi _{\sf eff}
\end{equation}
\blu{Information bounds:}\hspace{.75em}
  $ {\cal J}^{-1} $ (model~${\cal Q}$) $\ge \;{\cal I}^{-1}$ (model~${\cal Q}_{\nu_0}$). \\
\blu{Classical adaptivity} (necessary condition): \hfill\mbox{\tiny Stein (1956)}
\begin{equation}
    {\cal J}^{-1} = {\cal I}^{-1} 
    \iff \Lambda \perp \partial_2 {\cal Q}
    \iff \psi _{\sf eff}= \psi _{\sf class}
\end{equation}
ICs exist iff ${\cal J} > 0 $ iff
  $a' \Lambda \notin \clin\partial_2 {\cal Q}$\enskip $ \forall a\in\R^k, a\ne0$.
  \hfill\mbox{\tiny Shen (1995), Rieder (2000)}\\
\blu{Bounded} ICs exist \hspace{.25em}iff $ a' \Lambda \notin$ \pur{
  $L_1$-closure} $\cl_1\lin \partial_2 {\cal Q}$\enskip $ \forall a\in\R^k, a\ne0$.
\end{frame}
\subsection{1.3 Semiparametric regression}
\begin{frame}\blu{1.3 Semiparametric Regression:}\quad
  $P _{\theta,\nu}=Q ^{(W _{\theta,\nu},Z)}=$ law of observations
  $(W _{\theta,\nu},Z)$, where $Z$ is some $k$-dim.\ covariate and
  $W _{\theta,\nu}$ are the responses.
\par  \smallskip
 Optimally bounded ICs at $(\theta,\nu)$ are of (sufficient) form
\vspace{-\smallskipamount}\begin{equation}\label{spa:eq:Hopt}
   \varrho= (A \Lambda-\xi-a)\,\min \Bigl\{ 1,
                   \frac{b}{|A \Lambda-\xi-a|} \Bigr\}
\vspace{-\smallskipamount}\end{equation}
for some $b\in (0,\infty)$, $A\in \R ^{k \times k}$, $\xi\in \clin \partial_2 {\cal Q}$,
and some $a$.
\par
If the joint law of $(W _{\theta,\nu},Z)$ is distorted (errors-in-variables): $a\in\R^k$.
\par If only the conditional laws $Q ^{W _{\theta,\nu}|Z=z}(dw)$ may be distorted and
the marginal $Q ^{Z}(dz)$ is kept ideal (error-free-variables), then $a:\R^k\to\R^k$,
such that $\Ew (\varrho|Z)=0$.
\par \smallskip \blu{Remark}\quad \sz
Actually, condition $\varrho \perp \partial_2 {\cal Q}$ (infinite-dim.) allows only
approximations of the optimal~$\rho$:
Assuming a CONS $g_1,g_2,\dots$ of $\clin \partial_2 {\cal Q}$, one can prove the
existence of $A_m\in \R ^{k \times k}$, $\xi_m\in \lin \{g_1,\dots,g_m\}$,
and $a_m$, such that the IC $\varrho_m$ of form (\ref{spa:eq:Hopt})---now
    $\perp \partial_2 {\cal Q}$ weakened to $\varrho_m\perp g_1,\dots,g_m$---tend
in $L_2(P)$ to the optimal $\varrho$ (not necessarily of form (\ref{spa:eq:Hopt})).
\hfill \mbox{\tiny Shen (1995), Ruckdeschel, Hable, Rieder (2010)}
\end{frame}
\begin{frame}\blu{Cox regression:}\quad Response variables
   $W _{\theta,\nu}=(T_{\theta,\nu} \land C, {\bf1}_{\{T_{\theta,\nu}\le C\}})$
from survival times $T_{\theta,\nu}$ and a bounded censoring time $C$; 
$T_{\theta,\nu}$  and $C$ stoch.\ independent given~$Z$.
\par
The cumulative hazard function of $T_{\theta,\nu}\:|Z$ assumed of form $e ^{\theta'Z}\nu$
for some $\theta\in\R^k$ and unknown, abs.~continuous baseline hazard function~$\nu$.
\par Then the parametric scores function $\Lambda$ at $(\theta,\nu)$ is
\vspace{-\medskipamount}\begin{equation}
     \Lambda \bigl((y,\delta),z \bigr) =
   \delta z -z \hskip3\ei e ^{\theta 'z}\, \nu(y)
\vspace{-\medskipamount}\end{equation}
and $\partial_2 {\cal Q}=B \,L_2(\nu)$ for the operator $B$ defined by
\vspace{-\medskipamount}\begin{equation}
   B \zeta \colon \bigl((y,\delta),z \bigr) \longmapsto
    \delta \,\zeta(y) -  e ^{\theta 'z}
    \int _{[\hskip6\ei0,y]} \zeta \,d \nu \,,\hspace{1.5em}
    \zeta\in L_2(\nu)
\vspace{-\smallskipamount}\end{equation}
The projection on $\cl \partial_2 {\cal Q}$ equals $\Pi_2=B (B^* B)^{-1}B^*$,
     $B^*$ the adjoint, where $(B^* B)^{-1}B^*(y)=\Ew(Z|Y=y, \delta=1)$.
\qquad Estimation of~$\theta$, since $\Pi_2(\Lambda)\ne0$,
is not adaptive w.r.t.~$\nu$.
\hfill \mbox{\tiny Bickel, Klaassen et al.~(1993), van der Vart (1998)}
\par \medskip \sz \blu {Remark} To the IC~$\varrho$ of form~(\ref{spa:eq:Hopt}),
a robust version of the Cox PLE is constructed, using
the order statistics to $T_{\theta,\nu} \land C$,
as an M-estimator with the random weights
   $\min \bigl\{ 1, \frac{b}{|A \Lambda-\xi-a|} \bigr\}$
evaluated at a starting estimate $(\tilde \theta, \tilde \nu)$,
and a weighted Breslow estimate of $\nu$ employing the same weights.
\hfill \mbox{\tiny Ruckdeschel, Hable, Rieder (2010)}
\end{frame}
\subsection{1.4 Mixture models}
\begin{frame}\blu{1.4 Exponential mixture models:}\quad
   $Q _{\theta,\nu}(dx)= \int M_{\theta,u}(dx)\,\nu(d u)$, \enskip each
   $M_{\theta,u}(dx)$ a pm with $\mu$-density 
   $f(x,\theta,u)= \exp \{u'T_\theta(x)+S _{\theta}(x)-b(\theta,u)\}$
   and distribution $\nu(du)$ of the incidental parameter.
Setting dot=$\partial/\partial\theta$,
\vspace{-\smallskipamount}\begingroup \mathsurround0em\arraycolsep0em
\begin{eqnarray} \hspace{-1.5em}& \Ds
    \Lambda(X, \theta,\nu)= \dot T(X,\theta)'\Ew(U|T) +
    \dot S(X,\theta)-\Ew( \dot b(\theta,U)|T) & \\ \hspace{-1.5em}  & \Ds
    \partial_2 {\cal Q}= \bigl\{\, w(X)\in L_2\bigm| \Ew w(X)=0\,,\enskip
           \mbox{$w(X) $ is $\sigma(T)$-measurable}\,\bigr\} & \\
    \hspace{-1em} & \Ds \Pi_2 \colon h(X) \longmapsto
    \Ew (h(X)|T(X,\theta)) - \Ew h(X) &
\end{eqnarray}\endgroup
\par \vspace{-\smallskipamount}
where $\partial_2 {\cal Q}=\clin \partial_2 {\cal Q}= \cl_1\lin \partial_2 {\cal Q}$.
\hfill \mbox{\tiny Bickel, Klaassen et al.~(1993))}
\par Optimally robust IC of necessary and sufficient form~(\ref{spa:eq:Hopt}):
      \hfill\mbox{\tiny Shen (1995), Fischer (2006)}
\vspace{-1.25\smallskipamount}\begin{displaymath}
   \varrho=\Gamma \min \{1, \frac{b}{|\Gamma|}\}\,, \qquad \Gamma=A \Lambda-\xi-a \,,
   \quad\Lambda=\Lambda(X,\theta,)
\vspace{-1.25\smallskipamount}\end{displaymath}
with $\xi\in L_2(T(X,\theta))$ and $a\in\R$ determined 
such that $\Ew(\varrho|T)=0$.
\par \vspace{.5\smallskipamount}%
\blu{Special case:}\quad $T(X,\theta)=T(X)$ and $S(X,\theta)=\theta'S(x)$.  Then
\vspace{-\smallskipamount}\begin{equation}
   \Lambda= S-\Ew( \dot b|T),\quad
   \bar{\Lambda}= S-\Ew(S|T)
\vspace{-\medskipamount}\end{equation}
The conditional density of $X$ on $T=t$ w.r.t.\ $\mu(dx|T=t)$ not depending
on $\nu$, $\Lambda$ and $\bar{\Lambda}$ do not depend on $\nu$: 
\blu{classical adaptivity}. More generally,
  $\xi\in L_2(T(X))$ and $a\in\R$ 
such that $\Ew(\varrho|T)=0$ do not depend on $\nu$:
\blu{robust adaptivity} (\S2).\hfill\mbox{\tiny Shen (1995), Fischer (2006)}
\par \null\vfill \end{frame}
\subsection{1.5 Finite dimensional case}
\begin{frame}\blu{1.4 Finite Dimensional Case:}\quad
In case $\nu\in H_\theta \subset$ some $\R^m$, differentiability (1) is assumed with
\vspace{-\smallskipamount}\begin{equation}
  \partial_2 {\cal Q}= \{\,b' \Delta\mid b\in\R^m\,\}
\vspace{-\smallskipamount}\end{equation}
for some nuisance scores $ \Delta \in L_2 ^m(Q _{\theta_0,\nu_0})$,
$\Ew \Delta=0$, ${\cal D}:=\Cov \Delta>0$. \\
Fisher information at $(\theta_0,\nu_0)$ for the full parameter $(\theta,\nu)$ is
\vspace{-\smallskipamount}\begin{equation}
  {\cal H}= \Cov \begin{pmatrix}\Lambda\cr \Delta\end{pmatrix}
  = \begin{pmatrix} {\cal I}& {\cal C}\cr {\cal C}'& {\cal D}
     \end{pmatrix} \,, \qquad {\cal C}=\Ew \Lambda \Delta'
\vspace{-\smallskipamount}\end{equation}
where $\det {\cal H}= \det {\cal D} \det {\cal J}$,
      ${\cal J}= {\cal I}- {\cal C} \,{\cal D}^{-1}{\cal C}' $, and
      $\Pi_2 \Lambda= {\cal C}\, {\cal D}^{-1} \Delta$. Then
\par \vspace{-\smallskipamount}\hfill \mbox{\tiny Neyman (1951): $C(\alpha)$-tests}
\vspace{-2\medskipamount}\begin{equation}
     \psi _{\sf eff}= {\cal J} ^{-1}( \Lambda-{\cal C} \,{\cal D}^{-1} \Delta)
     \vspace{-.33\smallskipamount}\end{equation}
and $\psi _{\sf eff}=$ first $k$ coordinates of
   $ {\cal H}^{-1}\bigl( {\Lambda \atop \Delta}\bigr) = \,
   \psi _{\sf class}^{\sf full}$ for the full parameter.
\par \smallskip
  \blu{Adaptivity} $\!\iff \!\!$
  \blu{$ {\cal C}=\Ew \Lambda \Delta'=0$}\enskip
(symmetric in main/nuisance parameter).
\par \null\vfill \end{frame}
\section{2. Robust Adaptivity}
\subsection{2.1 Minmax MSE}
\begin{frame}\null
\blu{\bf Minmax MSE problems} \enskip
for AL estimators in robust neighborhood models:
   $$  \MSE_*(\psi,r)=\Ew |\psi|^2 + r^2 \, \omega_*^2(\psi) =\min {} ! $$
that, in addition to the asymptotic variance $\Cov \psi$, involve the
maximum asymptotic biasses $\omega_*$ generated by 
shrinking $r/\!\sqrt{n}\,$-neighborhoods about~$Q _{\theta_0,\nu_0}$,
  \blu{$\omega_* =$ sup-norm and variants},
  \hfill \mbox{\tiny e.g., integral of sectionwise supnorms, Rieder~(1994)}\\
and refer to the following two sets of ICs, respectively:
\par \medskip
1. model~${\cal Q} _{\nu_0}$ (no nuisance~$\nu$):
   $ \psi \in L_2^k(Q _{\theta_0,\nu_0})$, $\Ew \psi=0$,
   $\Ew \psi \Lambda'= \EM_k $\\[.5ex]
2. model~${\cal Q}$ (with nuisance~$\nu$): in addition, \enskip
   \blu{$ \Ew \psi g=0 $\enskip $\forall\,g\in \partial_2 {\cal Q}$}
\par \smallskip
Due to strict convexity, the minimizers $\varrho _1$ and $\varrho _2$,
respectively, are unique, and $ \mbox{minMSE1}\le \mbox{minMSE2}$.
\blu{Robust adaptivity} (extending classical):
\begin{equation}
   \mbox{minMSE1 = minMSE2}\iff \varrho_1=\varrho_2\iff
   \varrho_1 \perp \partial_2 {\cal Q}
\end{equation}
\blu{Nonadaptivity} (quantitative):
\qquad $\Ds\frac{\mbox{minMSE2}}{\mbox{minMSE1}}-1 $.
\end{frame}
\subsection{2.2 Symmetric location}
\begin{frame}\null \blu{\bf 2.2 Symmetric Location}
\hfill \mbox{\tiny Beran~(1974), Stone~(1975)}
\begin{equation}
  Q_ {\theta,f}(dx)=f(x-\theta)\,\lambda(dx),\quad \theta\in\R
\end{equation}
   $f$ symmetric,\: 
   ${\cal I}_f ^{\sf loc}=\int (\Lambda_f ^{\sf loc})^2\,f\,d\lambda<\infty$,\:
   $ \Lambda_f ^{\sf loc}=-{\dot f}\!/f$,\: $dF=f\,d\lambda$.\\
For $\theta_0=0$, $f=f_0$ fixed, \quad
    $\partial_2 {\cal Q}= \{\,g\in L_2(F)\mid
      \Ew g=0,\enskip \mbox{$g$ symmetric}\,\}$.\\
By symmetry,
  $ \Lambda_f ^{\sf loc}=- \dot f \!/f$ (odd) $\perp g $ (symmetric) 
  in~$L_2(F)$:\quad $\so $\\ 
  \qquad $ \Lambda_f ^{\sf loc}\notin \cl_1\lin\partial_2 {\cal Q}$\quad  
  and \enskip \blu{classical adaptivity} \enskip holds. 
\par \smallskip 
Robust ICs, for known $\nu_0=f$,
\hfill \mbox{\tiny Huber (1981), Hampel et al.\ (1985), Rieder (1994)}
   $$ \varrho(x)= A\,\Lambda_f ^{\sf loc}(x)\min \{1, c\,|\Lambda_f ^{\sf loc}(x)| ^{-1}\}$$
are all odd (like $\Lambda_f ^{\sf loc}$), hence $ \varrho \perp \partial_2 {\cal Q}$:
  $\so $ \blu{robust adaptivity} 
\par \medskip {\sz \blu{Remark}\quad
Adaptive constructions that not only achieve asymptotic linearity with the robust IC
in the ideal model but uniform asymptotic normality over shrinking neighborhoods
not yet solved completely. \hfill \mbox{\tiny Shen (1994), Stabla (2005)}}
\end{frame}
\subsection{2.3 Regression}
\begin{frame}
\null \blu{\bf 2.3 Regression and Scale}\hfill \mbox{\tiny Kohl (2005)} 
\begin{equation}
     Q _{\theta, \sigma}(dx,dy)=
     \frac{1}{\sigma}f\Bigl(\frac{y-x'\theta}{\sigma}\Bigr)\,\lambda(dy)K(dx)
\end{equation}
Assumptions:  \blu{$F$ symmetric,}
finite Fisher information of location ${\cal I}_f ^{\sf loc}$ and scale
     ${\cal I}_f ^{\sf sc}=\int (\Lambda_f ^{\sf sc})^2 \,dF $, where
     $ \Lambda_f ^{\sf sc}(u)=u \Lambda ^{\sf loc}(u)-1$;
     $ \int xx'\,K(dx)>0$. 
\par \smallskip
\blu{Classical adaptivity} holds 
(i.e., w.r.t.~$\sigma$ and w.r.t.~$\theta$)---due to symmetry of~$F$---and
extends to \blu{robust adaptivity} w.r.t.~$\sigma$ and w.r.t.~$\theta$,
in connection with the bias terms
\begingroup \mathsurround0em\arraycolsep0em
\begin{eqnarray}
    \omega_ {c,0}(\psi) &{}={}& \omega _{c,1}(\psi)= \Nsup _{x,u}|\psi(x,u)| \\
    \omega _{c,2}^2(\psi) &{}={}& \int \Nsup_u |\psi(x,u)|^2 \,K(dx)
\end{eqnarray}\endgroup
{\sz These biasses are generated by contamination
neighborhoods (Tukey, \blu{$*=c$}), which are unconditional (\blu{$t=0$}),
or errors-in-variables, or are average conditional,
error-free-variables, (\blu{$t=\alpha=1$}), respectively
by average square conditional neighborhoods (\blu{$t=\alpha=2$, $*=c$}).}
\hfill \mbox{\tiny Bickel (1984), Rieder (1987)}
\end{frame}
\begin{frame}
\null \blu{\bf Robust ICs for regression and scale}\quad {\fz
  $F$ symmetric, \quad $t=0$ and $t=\alpha=1$
\par \smallskip
$\theta$ main, $\sigma$ nuisance:
\begingroup \mathsurround0em\arraycolsep0em \vspace{-\medskipamount}
\begin{eqnarray}
   \varrho _{\sf rg}(x,u) &{}= {}& A_{\sf rg}x \Lambda_f ^{loc}(u) \,w _{\sf rg}(x,u)\\
   w _{\sf rg}(x,u) &{}= {}&
   \min \{\,1, b _{\sf rg}\,|A_{\sf rg}x \Lambda_f ^{\sf loc}(u)|^{-1}\} \\
   A_{\sf rg}^{-1} &{}= {}& \Ew xx' \Lambda_f ^{\sf loc}(u)^2 w _{\sf rg}(x,u)\\
   r^2 b_{\sf rg} &{}= {}& \Ew \bigl(|A_{\sf rg}x \Lambda_f ^{\sf loc}(u)|-b_{\sf rg}\bigr)_+
\end{eqnarray}\endgroup
$\sigma$ main, $\theta$ nuisance:
\begingroup \mathsurround0em\arraycolsep0em \vspace{-\medskipamount}
\begin{eqnarray}
  \varrho _{\sf sc}(u) &{}= {}& A_{\sf sc} (\Lambda_f ^{\sf sc}(u)-z_{\sf sc})w_{\sf sc}(u)\\
  w_{\sf sc}(u)&{}= {}& \min \{\,1, c_{\sf sc}\,|\Lambda_f ^{\sf sc}(u)-z_{\sf sc}|^{-1}\}
  \hspace{-1em}\\
                z_{\sf sc}   &{}= {}& \Ew \Lambda_f ^{\sf sc}w_{\sf sc}/\Ew w_{\sf sc} \\
            A_{\sf sc} ^{-1}   &{}= {}& \Ew (\Lambda_f ^{\sf sc}-z_{\sf sc})^2 w_{\sf sc}\\
       r^2 c_{\sf sc}  &{}= {}& \Ew \bigl( |\Lambda_f ^{\sf sc}-z_{\sf sc}|-c_{\sf sc}\bigr)_+
\end{eqnarray}\endgroup
Full parameter $(\theta,\sigma)$:\hspace{3.8em}
  $\varrho= \bigl({\varrho _{\sf rg}\atop \varrho _{\sf sc}}\bigr)$, \qquad
but weights $w _{\sf rg}$, $w_{\sf sc}$ both replaced by
\begingroup \mathsurround0em\arraycolsep0em
\begin{eqnarray}
  w(x,u)&{}={}& \min \bigl\{1, b \:\bigl| |A_{\sf rg}x|^2(\Lambda_f ^{loc}(u))^2 +
                        A_{\sf sc}^2 (\Lambda_f ^{\sf sc}(u)-z_{\sf sc})^2\bigr| ^{-1/2}
                        \bigr\} \\
  \noalign{\noindent where \vspace{-1\medskipamount}\nopagebreak}
  r^2b &{}={}&   \Ew \bigl( \:\bigl| |A_{\sf rg}x|^2(\Lambda_f ^{loc}(u))^2 +
                        A_{\sf sc}^2 (\Lambda_f ^{\sf sc}(u)-z_{\sf sc})^2\bigr| ^{1/2}-b\bigr)_+
\end{eqnarray}\endgroup   }{\fz
Especially, if $F={\cal N}(0,1)$,
then \pur{$|\varrho _{\sf rg}(x,u)|\sim 1/u $} \enskip ($x$ fixed, $|u|\to \infty$).}
\end{frame}
\begin{frame}\null\blu{\bf Regression with intercept as nuisance parameter}
\begin{equation}
   Q _{\theta,\,\mu}(dx,dy)= f(y-\mu-x'\theta)\,\lambda(dy)\,K(dx)
\end{equation}
Assumptions:  $F$ symmetric, ${\cal I}_f ^{\sf loc}<\infty$,
     $ \int xx'\,K(dx)>0$, $\int x\,K(dx)=0$.
\par \medskip
\blu{Classical adaptivity}, even if $K$ is asymmetric.
\par \smallskip
\blu{Robust adaptivity} for average square conditional neighborhoods\\
\qquad \blu{$t=\alpha=2$, $*=c$}, even if $K$ is asymmetric.
\par \smallskip
\blu{Robust adaptivity} for unconditional neighborhoods \enskip \blu{$*=c$, $t=0$}
and average conditional neighborhoods \enskip \blu{$*=c$, $t=\alpha=1$},
if \blu{$K$ symmetric.}
\par \smallskip
For asymmetric $K$, $*=c$, $t=0$ or $t=\alpha=1$, no robust adaptivity, since
\begin{equation}
   \Ew \varrho _{\sf rg} \,\Lambda_f ^{\sf loc} =
   A _{\sf rg} \Ew x \Lambda_f ^{\sf loc}(u)^2
   \min \{\,1, b _{\sf rg}\,|A _{\sf rg}\Lambda_f ^{\sf loc}(u)| ^{-1}\} \ne0
\end{equation}
For a 2-point asymmetric $K$, \blu{nonadaptivity} may be up to \blu{300\%}
\hfill \mbox{\tiny Kohl (2005)} 
\par \medskip 
{\fz
Robust ICs in model~${\cal Q}$ are of form (15), (16), (18) with $A_{\sf rg}x $
replaced by $A_{\sf rg}x + A_\mu$, and (17) by
\begingroup \mathsurround0em\arraycolsep0em \vspace{-\medskipamount}
\begin{eqnarray}
  A_{\sf rg} \Ew xx' \Lambda_f ^{\sf loc}(u)^2 w & {}= {} &
  \EM_k- A_\mu\Ew x' \Lambda_f ^{\sf loc}(u)^2 w \\
  A_\mu \Ew \Lambda_f ^{\sf loc}(u)^2 w & {}= {} &
  - A_{\sf rg} \Ew x \Lambda_f ^{\sf loc}(u)^2 w
\end{eqnarray}\endgroup }
\end{frame}
\subsection{2.4 ARMA}
\begin{frame}\null\blu{\bf\boldmath 2.3 ARMA$(p,q)$:\qquad
     $ \phi(B)(X_t- \mu)= \xi(B)V_t $}\quad $ t\in\Z$, \hfill $B$ backshift
\par \smallskip Innovations $V_t$ i.i.d.$\sim F$,
\quad ${\cal I}_F^{\rm loc}<\infty$, $\int u\,F(du)=0$, $ \int u^2 \,F(du)<\infty$.
\par 
Stationarity and invertibility assumption:
$\phi(z)\xi(z)\neq \;\forall\,|z|\le 1$, \\
$\phi$, $\xi$ relatively prime ($\Rightarrow$ positive Fisher information),
\enskip $\phi_p\,\xi_q\ne0$.
\par 
\blu{Influence} $\psi(x _{\le t})$ of observation~$x_t$
\blu{given the past} $x _{<t}:=(x_{t-1},x_{t-2},\dots)$.\\
Influence curves $\psi(x _{\le t})$ of AL estimators as in~(2),
but \blu{$\Ew \bigl(\psi(x _{\le t})\big|x_{<t}\bigr)=0$}
(stationary, ergodic martingale differences).
\hfill \mbox{\tiny Jeganathan (1982), Staab (1984), Rieder (2003)}\\
Differentiability~(1) now refers to transition densities of the 
ideal model~${\cal P}$.
\par \smallskip Joint law of $x _{\le n}$\,:\qquad
$    Q^{(n)}(dx _{\le n})= \prod\nolimits_{j=1}^n
     Q^{(n,j|<j)}(dx_j|x _{<j})Q^{(n,0)}(dx _{\le0})
$
\par \smallskip
\blu{Neighborhoods} ($*=c$, $t=\varepsilon$) of radius $r_n=r\,s_n= r\,n ^{-1/2}$
about the ideal \blu{transition distributions} $P^{(n,j|<j)}(dx_j|x _{<j})$
with \blu{contamination curve~$\varepsilon$:}
\begingroup \mathsurround0em\arraycolsep0em
\begin{eqnarray}\lefteqn{\hskip-2em
  Q^{(n,j|<j)}(dx_j|x _{<j}) = \null } \\
  && (1- r_n \:\varepsilon (x _{<j}))  P^{(n,j|<j)}(dx_j|x _{<j})
  + r_n \:\varepsilon(x _{<j}) M^{(n,j|<j)}(dx_j|x _{<j}) \nonumber
\end{eqnarray}\endgroup
where $M^{(n,j|<j)}(dx_j|x _{<j})$ any kernel,
initial distribution (of $x _{\le0}$) left ideal.
\par \quad \blu{$\alpha=1$:\enskip $\Ew \varepsilon\le1$,
\qquad $\alpha=2$:\enskip $\Ew \varepsilon^2\le1$}
\hfill \mbox{\tiny Bickel (1984), Rieder (1987) for regression}
\end{frame}
\begin{frame}
\blu{Bias terms} for $*=c$ and $t=\varepsilon$, respectively $t=\alpha=1,2$:
\begin{equation}
    \omega _{c,\varepsilon}(\psi) =
    \Ew \varepsilon(x _{\le0}) \sup \nolimits _{x_1}|\psi(x_1,x _{\le0})|
\end{equation}\begin{equation}
    \omega _{c,1}(\psi)= \Vert \psi \Vert _{\infty}\,,\quad
    \omega _{c,2}^2(\psi)= \Ew \sup \nolimits_{x_1}|\psi(x_1,x _{\le0})|^2
\end{equation}
\par
\blu{Transition scores:} \hfill
   $ \Lambda_1 = \Lambda_f ^{\sf loc}(V_1)\bigl(H_1',\tau\bigr)' $
\quad where $\tau=\phi(1)/\xi(1)$ \quad  and
\begin{equation}
 H_1' = \bigl(- B\phi^{-1}(B),\ldots,- B^p\phi^{-1}(B);
                     B\xi^{-1}(B),\ldots, B^q\xi^{-1}(B)\bigr)V_1
\end{equation}
Denoting \enskip ${\cal K}=  \Cov H_1$, \enskip Fisher information
is: \quad $ {\cal I}= {\cal I}_{F}^{\rm loc} \:\blu{\diag({\cal K}, \tau^2) }$
\par \smallskip
\qquad $\so$ \blu{classical adaptivity} (w.r.t.~$\mu$ and w.r.t.~$(\phi,\xi)$)
\par \medskip
Analogy to regression with intercept on identifying $H_1$ as regressor.
Robust ICs are of regression type form (15), (16), (18), (28), (29). \\
In the model with parameter $(\phi,\xi,\mu=0)$:
\begin{equation}
 \varrho _{c,\alpha}= AH_1 (\Lambda_f ^{\sf loc}(V_1)- \vartheta_\alpha)\,w_\alpha \,,
 \quad w_\alpha= \min \bigl\{\, 1,
\frac{\beta _{\alpha}}{|\Lambda_f ^{\sf loc}(V_1)- \vartheta_\alpha|}\,\bigr\}
\end{equation}
\par
\blu{$\alpha=1$:} \hspace{3em} $\beta_1=b/|AH_1|$, \hspace{.75em}
              $\vartheta_1=\vartheta_1(H_1)$ \hfill  \blu{Hampel-type}
\par \smallskip
\blu{$\alpha=2$:} \hspace{3em} $\beta_2=$ constant, \hspace{.5em}
              $\vartheta_2=$ constant \hfill \blu{Huber-type}
\end{frame}
\begin{frame}
\null \blu{\bf Robust Adaptivity for ARMA}
\par \medskip
1) Estimation of $(\phi,\xi)$, nuisance parameter~$\mu$:
\par \smallskip
\blu{Robust adaptivity} in case $\alpha=2$,
     in case $\alpha=1$ if~$F$ is symmetric.
\par \smallskip
In fact, $\Ew H_1=0$, 
and $H_1$,~$\Lambda_f ^{\sf loc}(V_1)$ are stochastically independent,
so 
\begingroup \mathsurround0em\arraycolsep0em
\vspace{-\smallskipamount}\begin{eqnarray*}\textstyle
  \Ew AH_1(\Lambda_f ^{\sf loc}- \vartheta_2)^2
      \min \bigl\{1, \frac{\beta_2}{|\Lambda_f^{\sf loc}-\vartheta_2|}\bigr\}
      &{}={}& 0       \qquad (\alpha=2)\\
\noalign{\vspace{-1 \smallskipamount}\noindent
         But \vspace{-\smallskipamount}\nopagebreak}
\textstyle \Ew AH_1(\Lambda_f^{\sf loc}-\vartheta_1)^2
      \min \bigl\{1, \frac{b/|AH|}{|\Lambda_f^{\sf loc}-\vartheta_1(H)|}\bigr\}
      &{}={}& 0       \qquad (\alpha=1)
\end{eqnarray*}\endgroup
where $ \vartheta_1(H_1)=\vartheta_1(-H_1)$, requires $\Lw_F(H_1)$, resp.~$F$,
to be symmetric.
\par \smallskip
\blu{Nonadaptivity} for AR(1), MA(1) with asymmetric~$F=\mbox{Gumbel}(\gamma,1)$,
  $\gamma=- \mbox{di}\Gamma(1)$ ($\Rightarrow\int v dF(v)=0$), at most 3\%!
\hfill \mbox{\tiny Kohl (2005)}
\par \bigskip
2) Estimation of~$\mu$ with nuisance parameter~$(\phi,\xi)$:
\par \smallskip
\blu{Robust adaptivity} for $\alpha=1,2$
\par \smallskip {\sz robust IC:  \qquad
   $\varrho _{c,12}= A \Lambda_f ^{\sf loc}(V_1) \min \{1,
   \beta _{12}\,|\Lambda_f ^{\sf loc}(V_1)|^{-1}\}$ \hfill $\alpha=1,2$}
\end{frame}
\subsection{2.5 ARCH}
\begin{frame}\null \blu{\bf 2.5 ARCH$(p)$: \quad
    $ X_t = \sigma(1+a_1 X_{t-1}^2+\ldots+a_p X_{t-p}^2)^{1/2}V_t $ }
    \qquad $ t\in\Z $
\par \medskip
Innovations $V_t$ i.i.d.$\sim F$, \enskip ${\cal I}_F^{\rm sc} < \infty$,
   $\int v dF(v)=0$, $\int v^2 dF(v)=1$\\
Stationarity, ergodicity: \qquad
  $\Ew\log V_t^2 + \log\sigma^2 + \log\max_j a_j< 0$
\par \medskip
Estimation of~$a$, nuisance parameter~$\sigma$:
\par \blu{No adaptivity}---neither classical nor robust ($*=c$, $\alpha=1$).
\par \smallskip
$\mbox{ARCH}(1)$ with $F=$ logNormal($\delta,\gamma$)
with $\delta = -e^{\gamma^2/2}$ ($\Rightarrow\int v dF(v)=0$):\\
\blu{Nonadaptivity} increases with~$r\in [\,0,\infty)$.
\hfill \mbox{\tiny Kohl (2005), MonteCarlo}\\
\qquad $a_1=1$, $\gamma=.5$: \enskip $.3 \uparrow .4$, \qquad
$a_1=10$, $\gamma=.5$: \enskip $25 \uparrow 160$. 
\par \bigskip \vfill 
\blu{\large\bf Conclusion}\quad {\sl
Classical adaptivity extends to robust adaptivity for neighborhoods
\enskip $*=c$, $\alpha=2$, \enskip for neighborhoods
\enskip $*=c$, $t=0$, $\alpha=1$ \enskip
some additional symmetry of the ideal model may be needed. \/}
\end{frame}
\section{3. Neighborhoods as Nuisance Parameter}
\subsection{3.1 Tangent balls}\begin{frame}
\blu{Neighborhood model}\quad
   $ {\cal Q}= \{Q\mid Q\in U_*(P_\theta,r),\:\theta\in \Theta\} $ of
neighborhoods about the elements of an ideal model
   ${\cal P}=\{P_\theta\mid \theta\in \Theta\}$.
Writing
\begin{equation}
  Q _{\theta,\nu}= P_\theta + \nu \,,\quad \blu{\nu:=Q-P_\theta}
  \quad \mbox{for \enskip $Q\in U_*(P_\theta,r)$}
\end{equation}
puts ${\cal Q}$ into semiparametric model form:
main parameter~$\theta$, nuisance parameter
  $\nu=Q-P_\theta\in H_\theta:=U_*(P_\theta,r)-P_\theta$;
  in particular, ${\cal P}={\cal Q}_{\nu_0=0}$.
\par \smallskip {\fz \blu{\bf Remark}
We assume a true~$\theta$ (idealistic approach), so the law~$Q$ may be referred
to this~$\theta$. Conversely, given~$Q$, the inclusion $Q\in U_*(P_\theta,r)$
may not define $\theta$ uniquely. }
\par \smallskip
\blu{Neigborhoods} $U_c(\theta,r)= \{(1-r)P+r\,M\mid \mbox{$M$ any probability}\}$
(convex contamination) and balls
   $U_*(\theta,r)= \{\,Q\mid d_*(Q,P_\theta)\le r\,\}$ 
in the Hellinger and total variation metrics, which are defined by
\begin{equation}
  \sqrt{2}\,d_h^2(Q,P)= \bigl\Vert \sqrt{dQ}\,-\sqrt{dP}\, \bigr\Vert_2 \,,\quad
  2\,d_v(Q,P)= \Vert dQ-dP \Vert_1
\end{equation}
\par \smallskip \blu{\bf 3.1 Proposition 1}\quad
{\sl
Fix $\theta_0$, $\nu_0=0$.
Then \enskip \blu{$\partial_2 {\cal Q}_* = r G_*$} for $*=h,v,c$,
where $G_*=$ all functions $g\in L_2(P _{\theta_0})$, $\Ew g=0$,
such that, respectively,
\begin{equation}
   (h)\enskip\Ew g^2\le 8 \qquad (v)\enskip\Ew|g|\le 2 \qquad
   (c)\enskip g\ge-1
\end{equation} \/}
In particular, if $r>0$:\hspace{.75em}%
       \pur{$\clin \partial_2 {\cal Q}_*=L_2(P _{\theta_0})$}, so
       $\Pi_2(\Lambda)=\Lambda$ and $\bar{\Lambda}=0$.
\end{frame}
\subsection{3.2 Semiparametric robust IC}\begin{frame}
We therefore dispense with the linear span and
define the \blu{sp-robust IC\/} 
\begin{equation}\label{nnn:eq:sp-robIC}
     \tilde{\varrho}_*:=  {\cal C} ^{-1} \tilde \Lambda \,,\qquad
     \tilde \Lambda=\Lambda- \widetilde{\Pi}_2(\Lambda)\,,\qquad
     {\cal C}= \Ew \tilde \Lambda \Lambda'
\end{equation}
in analogy to $\psi_{\sf eff}$, but employing nonlinear projection
\pur{$\widetilde{\Pi}_2 \colon L_2^k\to (r G_*)^k$} onto closed convex sets;
the radius $r$ is assumed so small that $\det{\cal C}\ne0$.
\par \medskip
{\sz \blu{\bf 3.2 Lemma 1}\label{nnn:lem:optapprox}
{\sl Let  $\tilde{\cal G}$, $\hat{{\cal G}}$, $\bar{{\cal G}}$ be some
nonempty closed convex: subset, cone, and linear subspace, respectively,
of some Hilbert space~${\cal H}$.
Then, for any $\kappa\in {\cal H}$, the unique best approximations
     $\tilde \kappa \in \tilde{\cal G}$, $\hat{\kappa}\in \hat{{\cal G}}$,
     $\bar{\kappa}\in \bar{{\cal G}}$ of~$\kappa$ are characterized by
\begin{equation}\label{nnn:eq:optapprox}
    \pur{ \langle \kappa-\tilde \kappa|g \rangle \le
    \langle \kappa-\tilde \kappa| \tilde \kappa\rangle },\hspace{1em}
    \langle \kappa-\hat \kappa|g \rangle \le
    \langle \kappa-\hat \kappa| \hat \kappa\rangle =0 \,,\hspace{1em}
    \langle \kappa-\bar \kappa|g \rangle =0
\end{equation}
for all $g$ in $\tilde{\cal G}$, $\hat{\cal G}$, and $\bar{\cal G}$, respectively. \/} }
\par \bigskip
\blu{\bf 3.2 Theorem 2} \quad
(h) \enskip {\sl If \enskip $ 8\,r^2<\min _{j=1,\ldots,k} {\cal I} _{j,j} $ \enskip then
    \enskip $ \tilde{\varrho}_h= {\cal I} ^{-1}\Lambda $\/}.
\par \smallskip
(v) \enskip {\sl If \enskip $2\,r< \min _{j=1,\ldots,k}\Ew|\Lambda_j|$ \enskip
then\/}\vspace{- \smallskipamount}
\begin{equation}\label{nnn:eq:spIC*v}
   \tilde \Lambda_j ^{(v)}= v_j' \lor \Lambda_j \land v_j '' \hspace{1.5em}
   \mbox{where}\hspace{.75em}
   \Ew (v_j'-\Lambda_j)_+=r=\Ew (\Lambda_j-v_j'')_+  \hspace{-1em}
\end{equation}
(c) \enskip {\sl If \enskip $r< -\max _{j=1,\ldots,k} \inf _{P _{\theta_0}}\Lambda_j$
\enskip then \vspace{- 1\smallskipamount}
\begin{equation}\label{nnn:eq:spIC*c}
    \tilde\Lambda_j ^{(c)}= ( \Lambda_j+r)\land u_j \qquad \mbox{where}\quad
    \Ew \bigl(( \Lambda_j+r)\land u_j \bigr)=0
\end{equation}   \/}
\end{frame}
\subsection{3.3 Comparison with robust IC}\begin{frame}\null
\blu{\bf Hellinger balls}\quad
Since $\MSE_h(\psi,r)=\tr \Cov \psi + 8 r^2 \maxev \Cov \psi$
and $\Cov \psi\ge {\cal I} ^{-1} \enskip
         \forall \, \mbox{ICs in model~${\cal P}$}$, the robust IC is
    $\varrho_h= {\cal I} ^{-1} \Lambda $ for all $ r\ge0$.
Thus \blu{$\tilde\varrho_h =\varrho_h$} whenever $\tilde\varrho_h $ is defined.
\par \smallskip {\fz
\blu{\bf Remark}\quad Despite $\varrho_h= \psi _{\sf class}$,
model~${\cal P}$ is not adaptive w.r.t.\ Hellinger neighborhoods since
   $ \MSE_h(\varrho_h,r)=\tr {\cal I} ^{-1} + 8 r^2 \maxev {\cal I} ^{-1}
     >\MSE_h(\varrho_h,0)$ for $r>0$. }
\par \bigskip \blu{\bf Contamination neighborhoods}\quad Risk
   $\MSE_c(\psi,r)= \Vert \psi \Vert_2^2+ r^2 \Vert \psi \Vert _{\infty}^2$
is uniquely minimized by the robust IC $\varrho_c$,
\begin{equation}
  \varrho_c=(A \Lambda-a) \min \bigl\{1, b\,|A \Lambda-a|^{-1}\bigr\}
\vspace{-1\medskipamount}\end{equation}
where\hfill \mbox{\tiny Rieder (1994)}\vspace{-1\medskipamount}%
\begin{equation}
   r^2 b= \Ew \bigl(|A \Lambda-a|-b \bigr)_+
\end{equation}
The sp-robust IC $\tilde \varrho_c$
(exchanging linear combination and clipping) has the coordinates 
\begin{equation}\hspace{-2em}
    \tilde \varrho_ {c,j}= C^-_{j,1}\, (\Lambda_1 + r) \land u_1
    + \cdots + C^-_{j,k}\, (\Lambda_k + r) \land u_k \hspace{-4em}
\end{equation}
with upper clipping constants $u_j$ from~(\ref{nnn:eq:spIC*c})
and $(C ^- _{j,i})={\cal C}^{-1}$ from (\ref{nnn:eq:sp-robIC}).
\par  \smallskip
In general, due to only one-sided (upper) bounds:
\quad \pur{$\MSE_c(\tilde \varrho_c,r)=\infty$\,!}
\end{frame}
\begin{frame}\null \blu{{\bf Total variation balls}---dimension $k=1$}
\par \smallskip
Robust IC minimizing
  $\MSE_v(\psi,r)= \Vert \psi \Vert_2^2+r^2 (\sup \psi -\inf \psi)^2$
is given by
\begin{equation}
   \varrho_v=c'\lor A\Lambda\land c''
\vspace{-\bigskipamount}\end{equation}
where \vspace{-\smallskipamount}
\begin{equation}
   r^2 (c''-c')=\Ew(c'-A\Lambda)_+=\Ew(A\Lambda-c'')_+
\end{equation}
\par \bigskip
\blu{\bf 3.3 Theorem 1}\quad {\sl
The sp-robust IC~$\tilde{\varrho}_v$ for $ r < \Ew \Lambda_+ $
coincides with the robust IC~$\varrho_v$ for
\vspace{-\smallskipamount}
\begin{equation}
  {\tilde r} = \sqrt{\frac{r}{v''_r-v'_r}}\,
\end{equation}
where
\begin{equation}
  \Ew (v_r'-\Lambda)_+ =  r =\Ew (\Lambda-v_r'')_+
\end{equation}\/}\hfill \mbox{\tiny Rieder (2000)}
\par 
{\sz
\blu{\bf 3.3 Example 2}\quad
In case $P_\theta={\cal N}(\theta,1)$, $\tilde\varrho_v$ turns out
pessimistic since
\par \medskip
\qquad  $\tilde r /r \ge 2.2 \enskip\forall r<1/\!\sqrt{2\pi}\,$,
\enskip  and $ \tilde r /r \uparrow \infty$\enskip
  as $r \downarrow 0$ or $\uparrow 1/\!\sqrt{2\pi}\,$  \\[1ex]
MSE-evaluation desirable. }
\end{frame}
\begin{frame}\null \blu{{\bf Total variation balls}---dimension $k>1$}
\par \smallskip
\par Robust IC minimizing
  $\MSE_v(\psi,r)= \Vert \psi \Vert_2^2+r^2 \omega _{v;s}^2$ for $s=2, \infty$
with\\[1ex]
\qquad $ \omega _{v;2}^2(\psi)= \sum _{j=1}^k (\sup \psi_j -\inf \psi_j)^2 $,
\enskip respectively\\[1ex] \qquad
  $ \omega _{v;\infty}^2(\psi)= \max _{j=1,\ldots,k}
                        (\sup \psi_j -\inf \psi_j)^2 $, \\[1.5ex]
has coordinates of form $ \varrho _{v,j}=c_j'\lor A_j\Lambda\land c_j''$
where, for variant $s=2$,\\[1ex]
\qquad
  $r^2 (c_j''-c_j')=\Ew(c_j'-A_j\Lambda)_+ =
            \Ew(A_j\Lambda-c_j'')_+ \quad \forall\,j=1,\ldots k$ \\[1ex]
respectively, for variant $s=\infty$,
\enskip $\forall\,j=1,\ldots k$,\\[1ex]
\qquad     $r^2 (c_j''-c_j')=
            \sum _{i=1}^k \Ew(c_i'-A_i\Lambda)_+
            = \sum _{i=1}^k \Ew(A_i\Lambda-c_i'')_+ $
\par \vspace{1.5\medskipamount}
Sp-robust IC $\tilde \varrho_v$ has the coordinates:
\begin{equation}
   \tilde{\varrho}_{v,j}= C^-_{j,1}\, v'_1 \lor \Lambda_1 \land v''_1
   + \cdots + C^-_{j,k}\, v'_k \lor \Lambda_k \land v''_k
\end{equation}
where $\Ew (v_j'-\Lambda_j)_+=r=\Ew (\Lambda_j-v_j'')_+$ and
      $(C ^- _{j,i})={\cal C}^{-1}$
      from (\ref{nnn:eq:sp-robIC}),~(\ref{nnn:eq:spIC*v}).
\par \smallskip
Thus the order of clipping and linear combination is interchanged again. \\
\blu{$\tilde \varrho_v$ is suboptimal but still sensibly robust.}
A MSE-comparison desirable.
\end{frame}
\subsection{3.4 A saddle point for testing convex sets}\begin{frame}\null
Dimension $k=1$.
Given any probability~$P$,
we consider local asymptotic alternatives along
tangents~$g\in L_2(P)$, $\int \!g \hspace{6\ei}dP= \Ew g = \langle g|1 \rangle=0$,
 \blu{
\begin{equation}
dP _{n,g} \approx (1+s_n g)\,dP \,,\qquad s_n=1/\!\sqrt{n}\,\hspace{-3em}
\end{equation}  }
{\fz E.g., by $P$-densities:
  $ \bigl(\,\frac{1}{2}sg+ (1-\frac{1}{4}s^2\Vert g \Vert^2) ^{1/2} \,\bigr)^2$,
  or simply: $1+sg$, if $\Vert g \Vert_\infty<\infty$. }
\enskip Observations $x_1,\ldots,x_n$ i.i.d.~$\sim P _{n,g}$.
\par \smallskip
Let $G_0,G_1\subset L_2\cap \{\Ew=0\}$, $G_0\cap G_1=\emptyset$.
Fix any $g:=(g_0,g_1)\in G_0 \times G_1$.
The simple asymptotic testing problem $H_{g_0}$ vs.~$K_{g_1}$
at level~$\alpha\in (0,1)$ is:
\begin{equation}
  \liminf_{n\to \infty} \int \!\delta_n \, dP_{n,g_1}^n = \max{!}
  \qquad \mbox{s.t.} \quad
  \limsup_{n\to \infty} \int \! \delta_n \, dP_{n,g_0}^n \le \alpha
\end{equation}
Denoting $g_{10}:=g_1-g_0$,
the \blu{optimal test} is $\delta_{g}=(\delta_{n,g})$,
\begin{equation}
   \delta_{n,g}=  {\bf 1}\biggl (s_n \smash{\sum_{i=1}^n}
   g_{10}(x_i) > \Vert g_{10}\Vert \hspace{6\ei} u_{\alpha}
   + \langle g_{10}|g_0\rangle  \biggr )
\end{equation}
{\fz where $\Vert \gP\Vert=\Vert \gP\Vert_2= \langle \gP|\gP \rangle ^{1/2}$,
and $u_{\alpha}=$ standard normal upper $\alpha$~point: $\Phi(-u_\alpha)=\alpha$.}\\
$\delta_{g}$~achieves asymptotic size~$\alpha$  
and power~$\Phi(-u_{\alpha}+\Vert g_{10}\Vert\,)$ under $H_{g_0}$, $K_{g_1}$.\\
The tests~$\delta_{n,g}$ are unique up to terms $\to 0$ in $P^n$-probability.
\end{frame}
\begin{frame}
The maxmin asymptotic testing problem
   $H_{G_0}$ vs.~$K_{G_1}$ at level $\alpha\in (0,1)$ is
\begin{equation}
      \inf_{g_1\in G_1}\liminf_{n\to \infty}
      \int \!\delta_n \, dP_{n,g_1}^n  = \max {!} \quad
      \mbox{s.t.} \enskip
   \sup_{g_0\in G_0}\limsup_{n\to \infty}
   \int \!\delta_n \, dP_{n,g_0}^n  \le  \alpha
\end{equation}
Assume now $G_0, G_1$ closed, convex.
Pass to \blu{$G _{10}:=G_1-G_0$}, which set is convex,
but need not be closed if $\dim L_2(P)>1$.
We assume \blu{$G _{10}$ closed} and pick \enskip $q _{10}:=q_1-q_0$ \enskip
the unique minimum norm element of $G _{10}$.
\par \smallskip
\blu{{\bf 3.4 Theorem 1}[\,saddle point for testing\,]}\quad {\sl \\
Then the maxmin asy.\ testing problem
  $H_{G_0}$ vs.~$K_{G_1}$ at level $\alpha$ has saddle point
  $(q,\delta_q)$, 
and the maxmin asy.~power $= \Phi(-u_{\alpha}+\Vert q_{10}\Vert\,)$.
\par
Any other pair $g=(g_0,g_1)$ in $G_0 \times G_1$
achieving $g _{10}=q _{10}$ 
also provides a saddle point $(g,\delta_g)$,
and necessarily  $\delta_g=\delta_q$. \/}
\par \smallskip {\sz
\blu{\bf Proof}\enskip Based on LAN, this is the \blu{statistical equivalent}
of the first characterization in (\ref{nnn:eq:optapprox})
with $\kappa=0$, for the minimum norm element of closed convex sets. }
\par \smallskip
Given some scores
  $\Lambda\in L_2 (P _{\theta_0})$, $\int \!\Lambda\,dP _{\theta_0}=0$,
and $\tau \in \R $, $\ne0$, enlarge the parametric alternatives \enskip
  $ dP _{\theta_0+s_n\tau} \approx (1+ s_n\tau \Lambda)dP _{\theta_0}$
to $P _{n,g}$, by the nuisance parameter $g\in G_0$, respectively
 \blu{$g\in \tau \Lambda + G_1$}.
Then 
\begin{equation}\pur{
  q _{10}= \tau\Lambda - \widetilde{\Pi}_2(\tau\Lambda | G_0-G_1)}
\end{equation}
\end{frame}
\subsection{3.5 Robust asymptotic tests}\begin{frame}\null
To test neighborhoods
  $U_*(\theta_0, s_n r_0)$ and $U_*(\theta_0+s_n\tau, s_n r_1)$
about $P=P _{\theta_0}$ and $P _{\theta_0+s_n\tau}$ of type $*=h,v,c $
with possibly different radii $s_n r_0$ and $s_n r_1$, respectively,
employ the tangent balls $G_*$ defined in (37) and put
\begin{equation}
  G _{*,0}=r_0 G_* \,, \qquad G _{*,1}= \tau \Lambda + r_1 G_*
\end{equation}
Abbreviate $H_*:=H_{G_{*,0}}$ and~$K_*:=K_{G_{*,1}}$.
\par \smallskip
\blu{{\bf 3.5 Theorem 1}\enskip [\,Hellinger balls, $*=h$\,]}\quad
{\sl
Let $ 8 \,r^2 < \tau^2\,\Vert \Lambda\Vert^2$. \enskip
Then the least favorable tangent pair
     $q_h=(q_{h,0},q_{h,1})$ in $G_{h,0}\times G_{h,1}$ is unique,
\begin{equation}
      q_{h,0}=r_0 \hspace{6\ei}\gamma \hspace{3\ei}\Lambda \,, \quad
      q_{h,1}   = \tau \hspace{3\ei}\Lambda -
      r_1 \hspace{6\ei}\gamma \hspace{3\ei}\Lambda \, \quad
      \mbox{where}\enskip \gamma= \sqrt{8}\:\Vert \Lambda\Vert^{-1}
\end{equation}
The maxmin 
test~$\delta_{q_h}=(\delta_{n,q_h})$ for~$H_h$ vs.~$K_h$ is
\begin{equation}
   \delta_{n,q_h} =
   {\bf 1}\biggl( s_n \Vert \Lambda \Vert ^{-1} \smash{\sum_{i=1}^n}
   \Lambda(x_i) >  u_{\alpha} + \sqrt{8}\,r_0 \biggr)
\end{equation}
\hfill Maxmin asymptotic power
  $ = \Phi \bigl(-u_{\alpha}+ \tau \hspace{6\ei}\Vert \Lambda\Vert
    - \sqrt{8}\,r\,\bigr) $. \/}
\par \smallskip {\fz \blu{\bf Remarks}\enskip 
a) Despite of classical test statistics, no adaptivity w.r.t.\ Hellinger balls. \\
b) No Huber--Strassen least favorable pairs $dQ_1=\pi\,dQ_0$ to compare with:
\hfill \mbox{\tiny Birg\`e (1980)}
\begin{equation}
      Q _{0}(\pi>t)\ge Q'(\pi>t),\; 
      Q _{1}(\pi>t)\le Q''(\pi>t)\quad
      \forall Q'\in {\cal Q}_0, Q''\in {\cal Q}_1 ,
      \forall\,t>0
\end{equation} }
\end{frame}
\begin{frame}\null\blu{{\bf 3.5 Theorem 2}\enskip
  [\,Total variation balls, $*=v$\,]}\quad
{\sl Let $2\,r < \tau \hspace{6\ei}\Ew |\Lambda|$.
\par 
{\sf a)}\enskip   Then a least favorable tangent pair
  $q_v=(q_{v,0},q_{v,1})$ in $ G_{v,0}\times G_{v,1}$ is given by
  \enskip $ q_{v,0} = r_0 \hspace{6\ei}\tilde{g}_v $,
          $ q_{v,1}= \tau \hspace{3\ei}\Lambda -r_1 \hspace{6\ei}\tilde{g}_v $
\enskip for the tangent $\tilde{g}_v$ defined by
\begin{equation}
    r\hspace{6\ei}\tilde{g}_v =
    \tau \hspace{6\ei}(\Lambda - v'' \hspace{6\ei})_+
    - \tau \hspace{6\ei}(v'- \Lambda)_+
\end{equation}
and $v'<0<v''$ determined by $
   \tau \hspace{3\ei}\Ew ( v' - \Lambda)_{+} = r
   = \tau \hspace{3\ei}\Ew ( \Lambda - v'' \hspace{6\ei})_{+}
   $. \\
A tangent pair $g _{v,0}=r_0g_0$, $g _{v,1}=\tau \Lambda-g_1$ is least favorable iff
\begin{equation}
    r_0 \hspace{6\ei} g_0^+ + r_1 \hspace{6\ei} g_1^+ =
    \tau \hspace{6\ei}(\Lambda - v'' \hspace{6\ei})_+  
    \,, \qquad r_0 \hspace{6\ei} g_0^- + r_1 \hspace{6\ei} g_1^- =
    \tau \hspace{6\ei}(v'- \Lambda)_+   
\end{equation}
With $ \Lambda^{(v)}: = v'\lor \Lambda\land v''$, the maxmin 
test~$\delta_{q_v}=(\delta_{n,q_v})$ for~$H_v$ vs.~$K_v$ is
\vspace{- \smallskipamount}\begin{equation}
   \delta_{n,q_v} =  {\bf 1}\biggl( s_n \smash{\sum_{i=1}^n}
   \Lambda^{(v)} (x_i) >  \Vert \Lambda^{(v)} \Vert
   \hspace{6\ei} u_{\alpha} + r_0 \hspace{6\ei}
   (v'' - v' \hspace{4\ei})   \biggr )
\vspace{-\bigskipamount}\end{equation}
\hfill Maxmin asy.~power
  $ = \Phi \bigl(-u_{\alpha}+ \tau \hspace{6\ei}\Vert \Lambda^{(v)}
    \Vert  \,\bigr) $.
\par \medskip {\sf b)}\enskip
$\delta_{q_v}$ coincides with the robust asy.~test based on
least favorable probability pairs
for~$U_v \bigl(P_{\theta_0};r_0/\!\sqrt{n}\,\bigr)$ vs.\
    $U_v \bigl(P_{\theta_0+\tau\!/\!\sqrt{n}\,};r_1/\!\sqrt{n}\,\bigr)$,
hence maximizes the asy.~minimum power over
  $U_v \bigl(P_{\theta_0+\tau\!/\!\sqrt{n}\,}; r_1/\!\sqrt{n}\,\bigr) $
subject to asy.~maximum size~$\le\alpha$
over~$U_v \bigl(P_{\theta_0};r_0/\!\sqrt{n}\,\bigr) $.    }
\end{frame}
\begin{frame}\null
\blu{{\bf 3.5 Theorem 3}\enskip [\,Contamination, $*=c$\,]}\quad
{\sl Let
   $ r_0 < \Ew \bigl(\tau \hspace{3\ei}\Lambda-(r_1-r_0)\bigr)_+ $.
\par {\sf a)}\enskip
The least favorable tangent pair~$q_c=(q_{c,0},q_{c,1})$
     in~$G_{c,0}\times G_{c,1}$ is unique,
\vspace{- \smallskipamount}\begin{equation}
\hspace{-2em}
   q_{c,0} = \tau\hspace{6\ei}
             (\Lambda- c'' \hspace{6\ei})_+ -r_0 \,, \quad
   q_{c,1} = 
   \tau \hspace{3\ei}\Lambda
     + \tau \hspace{6\ei} (c'-\Lambda)_+ -r_1
\hspace{-4em}\end{equation}
where $c'<z:=(r_1-r_0)/\tau < c''$ are determined by
\hfill $\Ew q_{c,0}=\Ew q_{c,1}=0$.
\par \smallskip
Based on $ \Lambda^{(c)}: =  c'\lor \Lambda\land c'' - z $, 
the maxmin 
test~$\delta_{q_c}=(\delta_{n,q_c})$ for~$H_c$ vs.~$K_c$ is
\vspace{-\medskipamount}\begin{equation}
   \delta_{n,q_c} =  {\bf 1}\biggl( s_n \smash{\sum_{i=1}^n}
   \Lambda^{(c)} (x_i) >  \Vert \Lambda^{(c)} \Vert
   \hspace{6\ei} u_{\alpha} + r_0 \hspace{6\ei} (c'' - z) \biggr)
\hspace{-2em}\vspace{-\medskipamount}\end{equation}
\hfill Maxmin asy.~power
  $ = \Phi \bigl(-u_{\alpha}+ \tau \hspace{6\ei}\Vert \Lambda^{(c)}\Vert
    \,\bigr) $.
\par \medskip {\sf b)}\enskip
$\delta_{q_c}$ coincides with the robust asy.~test
based on least favorable probability pairs
for~$U_c \bigl(P_{\theta_0};r_0/\!\sqrt{n}\,\bigr)$ vs.\/~{\rm
    $U_c \bigl(P_{\theta_0+\tau\!/\!\sqrt{n}\,};r_1/\!\sqrt{n}\,\bigr)$,}
hence maximizes the asy.~minimum power over
  $U_c \bigl(P_{\theta_0+\tau\!/\!\sqrt{n}\,}; r_1/\!\sqrt{n}\,\bigr) $
subject to asy.~maximum size~$\le\alpha$
over~$U_c \bigl(P_{\theta_0};r_0/\!\sqrt{n}\,\bigr) $. \/}
\par \vfill
\hfill \mbox{\tiny Huber (1964), (1968), Huber--Carol (1970),
                   Huber--Strassen (1973), Rieder (1978), (2000)}
\end{frame}
\begin{frame}\null \blu{\large\bf Summary}
\par \bigskip
ESTIMATION
\par \medskip
Hellinger ($*=h$): SpM (semiparametric method) yields the optimally robust IC.
\par \smallskip
Total variation ($*=v$), parameter dim $k=1$:
SpM yields a suboptimal IC of optimally robust form (for a different radius).
\par \smallskip
Total variation ($*=v$), parameter dim $k>1$:
SpM eases the problem by exchanging the order of clipping and linear
combination of coordinates.
The sp-robust IC thus obtained is 
reasonably robust under MSE. 
\par \smallskip
Contamination ($*=c$): SpM fails, yielding unbounded ICs, ${\rm MSE}=\infty$.
\par \bigskip
TESTING \quad a one-dimensional parameter
\par \medskip
Total variation, contamination ($*=c,v$): SpM yields the optimally
robust---maxmin---asymptotic tests of {\sl Huber--Strassen\/} form.
\par \smallskip
Hellinger ($*=h$): SpM yields a maxmin asymptotic test---{\sz although,
at finite sample size, no Huber--Strassen pairs exist}.
\par \vfill\end{frame}
\section{4. Uniform Asymptotic Normality}
\subsection{4.1 Adaptive estimators}\begin{frame}
Adaptive constructions by Beran (1976), Krei\ss\ (1987) for ARMA, and
by Drost, Klaassen, Wercker (1997, 1998) for ARCH, GARCH, TAR,
such that for all $F$, $\theta$
\begin{equation}
  \Lw _{F,\theta} \bigl\{ (n\,{\cal I} _{F,\theta}) ^{1/2}
    (S_n-\theta)\bigr\} \longrightarrow {\cal N}(0, \EM_k)
\end{equation}
Adaptation w.r.t.\ symmetric innovation distribution.
\par \smallskip
{\bf \textcolor{purple}{Nonuniformity}} \hfill \mbox{\tiny Klaassen~(1980)}\\
{\sl 1-dim location, $F$ symmetric, ${\cal I} ^{\sf loc}_F<\infty$,
$S_n \colon \R^n\to\R$ translation equivariant, sample size~$n$ fixed.
Then \enskip $\forall\,\varepsilon>0$ $\forall\,x>0$
\begin{equation}\label{Klaassen}
    \inf _{G\in B_c ^{\sf s,i}(F,\varepsilon)}
    G^n \bigl\{|(n\,{\cal I}_G) ^{1/2}S_n|\le x\bigr\} =0 < 2\,\Phi(x)-1
\end{equation}
where}\qquad $   B_c ^{\sf s,i}(F,\varepsilon)=
         \bigl\{(1-\varepsilon)F+\varepsilon H \bigm|
         \mbox{$H$ symmetric},\; {\cal I} ^{\sf loc}_H<\infty \bigr\}
                $
\par \bigskip
Extensions to other models? Practical use of adaptive estimators?
Robustness? \hfill \mbox{\tiny Bickel (1981), (1982), Huber (1996)}
\vfill
\end{frame}
\subsection{4.2 Models, Fisher information}\begin{frame}\null
\blu{\bf Models}\quad
Location, scale (nonidentifiable), linear regression, ARMA\\
having a finite Fisher information of the form
\begin{equation}\blu{
     {\cal I} _{F,\theta} = {\cal I}_F ^{\sf loc/sc}
     \sigma_F^2 \:{\cal K} _{\theta} }
\end{equation}
{\fz Factor $\sigma_F^2 = \int v^2 \,F(dv)$, where $\mu_F=\int v \,F(dv)=0$, appears only
in MA, AR, ARMA.}
\par \mbox{\tiny Huber (1981)}
\begin{equation}\blu{
     {\cal I}_F ^{\sf loc}:}= \Nsup _{\varphi\in {\cal C}^1_c}
       \bigl(\Tint {\dot \varphi}\,dF\bigr)^2 \!\big/\!\!\Tint \varphi^2 dF 
\smallskip \end{equation} 
{\fz ${\cal C}^1_c:=$ all continuosly differentiable functions of compact support.
Then: ${\cal I}_F ^{\sf loc}<\infty$ iff $dF=f \,d\lambda$, $f$ abs.~continuous
and  $\int (\Lambda_F ^{\sf loc})^2 \,dF<\infty$, in which case
     $ {\cal I}_F ^{\sf loc} = \int (\Lambda_F ^{\sf loc})^2 \,dF$.}
\par \mbox{\tiny Ruckdeschel, Rieder (2010)}
\begin{equation}\blu{
     {\cal I}_F ^{\sf sc}:}= \Nsup _{\varphi\in {\cal C}_{1,c}}
       \bigl(\Tint v{\dot \varphi}(v)\,dF\bigr)^2 \!\big/\!\!\Tint \varphi^2 dF 
     \hspace{-3em} 
\end{equation} 
{\fz
  ${\cal C}_{1,c}:=$ all functions with continuous derivative of compact support.
  Then:
  ${\cal I}_F ^{\sf sc}<\infty$ iff $dF=f \,d\lambda$ on $\R \setminus \{0\}$,
  $v \mapsto v \,f(v)$ is abs.~continuous and
  $\int _{\ne0} (\Lambda_F ^{\sf sc})^2 \,dF<\infty$ where
  $\Lambda_F ^{\sf sc}= v \,\Lambda_F ^{\sf loc}-1$,
  in which case
  $ {\cal I}_F ^{\sf sc} = \int _{\ne 0} (\Lambda_F ^{\sf sc})^2 \,dF$. }
\par \medskip
  $\so$ \qquad \pur{ ${\cal I}_F ^{\sf loc/sc}$ is 
  convex and weakly l.s.c.\, but not u.s.c.\hskip24\ei ! } 
\end{frame}
\subsection{4.3 Lower bounds in Kolmogorov metric}\begin{frame}
Kolmogorov metric = sup-norm distance between c.d.f.'s on~$\R^k$
\par \smallskip
\blu{\bf 4.3 Theorem 1} (location, scale, linear regression, MA)\\ {\sl
Assume ${\cal I} ^{\sf loc/sc}_F<\infty$,
$S_n \colon \R^n\to\R^k$ any estimator, $n$ fixed.
Then $\:\forall\,\varepsilon>0$\vspace{-\smallskipamount}
\begin{equation}
  \sup_{G\in B_c ^{\sf s,i}(F,\varepsilon)} d_\kappa  \Bigl( \Lw _{G,\theta} \bigl\{
   (n\,{\cal I}_{G,\theta}^{1/2}(S_n-\theta)\bigr\}, {\cal N}(0,\EM_k)\Bigr)
   \ge 1 - \frac{1}{2^k} - \kappa_n
\end{equation}
where \vspace{-\medskipamount}
\begin{equation}
   \kappa_n:=d_\kappa \Bigl( \Lw _{F,\theta} \bigl\{
   (n\,{\cal I}_{F,\theta}^{1/2}(S_n-\theta)\bigr\}, {\cal N}(0,\EM_k)\Bigr)
\end{equation}
and $ B_c ^{\sf s,i}(F,\varepsilon)= $ all $(1-\varepsilon)F+\varepsilon H$
with $H$ symmetric, ${\cal I} ^{\sf loc/sc}_H<\infty$, and,
       in case MA, in addition $\mu_H=0$, $\sigma_H^2\in (0,\infty)$. \/}
\par \medskip {\sz
\blu{\bf Remarks} a) No equivariance, no symmetry assumptions.
\par \smallskip
b) Use $G_m=(1-\varepsilon_m)F + \varepsilon_m/2 \,\bigl(
         {\cal N}(-a,\sigma_m^2)+ {\cal N}(a,\sigma_m^2)\bigr)$ with
       $\varepsilon_m, \sigma_m^2\to0$ such that
   ${\cal I}^{\sf loc/sc} _{G_m}\to \infty$ and, in case of MA,
   $\sigma^2 _{G_m}\to \sigma^2 _{F}$. In these models, the
joint law of observations is $d_v$-continuous in the innovation distribution.
Pass to $d_\kappa$, which is scale invariant and metrizes weak convergence
to ${\cal N}(0,\EM_k)$.
\par \smallskip
c) \hspace{9.75em}\mbox{
   $1-2 ^{-k}= d_\kappa \bigl(1_0, {\cal N}(0,\EM_k)\bigr)$ }  }
\end{frame}\begin{frame}\null
\blu{\bf 4.3 Theorem 2} (location, scale, linear regression, MA, AR, ARMA)
\par \smallskip  {\sl
If \quad $ \Lw _{F,\theta} \bigl\{(n\,{\cal I}_{F,\theta}^{1/2}(S_n-\theta)\bigr\}
     \longrightarrow {\cal N}(0,\EM_k) $ \quad then, for any $\varepsilon_n \to0$,
\begin{equation}
   \liminf _{n\to \infty} \sup_{G\in B_c ^{\sf s,i}(F,\varepsilon_n)}
   d_\kappa  \Bigl( \Lw _{G,\theta} \bigl\{(n\,{\cal I}_{G,\theta}^{1/2}
         (S_n-\theta)\bigr\}, {\cal N}(0,\EM_k)\Bigr)
    \ge 1 - \frac{1}{2^k}
\end{equation}
where $ B_c ^{\sf s,i}(F,\varepsilon)= $ all $(1-\varepsilon)F+\varepsilon H$
with $H$ symmetric, ${\cal I} ^{\sf loc/sc}_H<\infty$, and,
       in the cases MA, AR, ARMA, in addition $\mu_H=0$, $\sigma_H^2\in (0,\infty)$.\\
In the case of AR, ARMA, the functions $S_n$ are required to be continuous. \/}
\par \medskip {\sz
\blu{\bf Remarks} a) Adaptive constructions~$S_n$ are smooth in the observations.
\par b) In AR, ARMA, i.e. MA($\infty$), the joint law of the observations is
not $d_v$-continuous in the innovation distribution. 
Instead, we derive bounds in~$L_2$ which translate to Prokhorov distance~$d_\pi$ via
\hfill \mbox{\tiny Strassen (1965)}
\begin{equation}
   d_\pi \bigl(\Lw(Y),\Lw(X)\bigr)\le \sqrt{{\Vert Y-X\Vert}_2}\,
\end{equation}
Invoke continuity of $S_n$ and, again, switch to~$d_\kappa$.
\par c) ARCH? GARCH?  }
\end{frame}
\subsection{4.4 Continuity of maximum risk}\begin{frame}
Let $({\cal M},d)$ be any metric space, balls $B(F,r)$ (open/closed). \\
For any given function $\alpha \colon {\cal M} \to\R$ consider
\begin{equation}
    \beta(F,r):= \sup \{\, \alpha(G)\mid G\in B(F,r)\,\}
\end{equation}
which, for fixed~$F$, increases in~$r$.
\par \smallskip
\blu{\bf 4.4 Lemma 1} {\sl The function~$\beta$ satisfies
\begin{equation}
  \beta(F,r-0)\le \liminf _{G\to F} \beta(G,r) \le
                  \limsup _{G\to F} \beta(G,r) \le \beta(F,r+0)
\end{equation}
with ``$=$'' except for countably many values of $r$, depending on~$F$. \/}
\par \smallskip
Follows from \quad $B(F,r-\delta)\subset B(G,r) \subset B(F, r+\delta)$ \enskip
if \enskip   $\delta=d(G,F)$.
\par \bigskip {\fz \blu{\bf Remarks}
a) Robust risk (max\,Var, max\,MSE, min\,FisherInfo) continuous.  \\
b) Weak dependence of robust estimators and minmaxrisk on
the unknown radius~$r$ of neighborhoods as a nuisance parameter.
\hfill \mbox{\tiny Rieder, Ruckdeschel, Kohl (2008)} \\
c) Based on uniform tightness of the empirical process, uniformly asymptotically
normal constructions of robust estimators in the independent case,
\begin{equation}
     \Lw _{Q_n^n} \bigl\{n ^{1/2}(S_n-T(Q_n))\bigr\}  \longrightarrow
     {\cal N}(0, \Cov \varrho_\theta)
\end{equation}
for all sequences $Q_n$ out of neighborhoods $U_*(\theta, r_n)$ about $P_\theta$,
$r_n=r\,n^{-1/2}$, $0<r<\infty$. \\
d) Difficulties under dependence; need neighborhoods smaller than (2.30).  }
\end{frame}
\section{5. One-Sided Inference on Tangent Cones} 
\subsection{5.1 Tangent Cones And Spaces}\begin{frame}
\blu{Functional} $T \colon {\cal P}\longrightarrow \R$,
defined on a family ${\cal P}$ of pm's on some sample space $(\Omega, {\cal A})$.
Observations $x_1,\dots,x_n$ i.i.d. $\sim$ any $P \in {\cal P}$.
\par
Want most accurate tests and confidence statements about $T(P)$.
\par
Fix any $P =P_0 \in {\cal P}$.
Local alternatives at $P$ within ${\cal P}$ are defined by
\begin{equation}\label{ccc:eq:Pgs}
\sqrt{\!\smash{dP _{g,s}}\vphantom{dP}}\, = \bigl(1+ \Tfrac{s}{2}g \bigr) \sqrt{\!dP}\, + \Lo(s)
\quad \mbox{as $s \downarrow 0$}
\end{equation}
\blu{Tangent set}\quad ${\cal G}$ of all $g\in L_2(P)$, $g \perp 1$, 
                  $P _{g,s}\in {\cal P}$ for small $s>0$.\\
  ${\cal G}$ is a cone, vertex at $0$
  {\fz (i.e., $\gamma g\in {\cal G}$ whenever $g\in {\cal G}$, $\gamma\ge0$)},
  and will be assumed also convex
  {\fz (i.e., $\gamma_1 g_1+\gamma_2 g_2\in {\cal G}$
  for all $g_0,g_1\in {\cal G}$, $\gamma_0, \gamma_1\ge0$)}.
\par \smallskip
\blu{Differentiability} of $T$: There is some $\kappa \in L_2(P)$ such that for all $g\in {\cal G}$,
\begin{equation}\label{ccc:eq:Tkappa}
   T(P _{g,s})=T(P)+s \langle \kappa|g \rangle+\Lo(s) \quad \mbox{as $s \downarrow 0$}
\end{equation}
$\kappa$ is nonunique, but $\bar{\kappa}=$ the orthoprojection of $\kappa$ onto $\clin{\cal G}$ is unique.\\
In addition, let $\hat \kappa=$ the (nonorthogonal) projection of $\kappa$ onto $\cl {\cal G}$.
\par \medskip
\pur{Literature} {\fz
The $*$-Theorem 25.20, for ${\cal G}$ a cone, 
and LAM-Theorem 25.21, for ${\cal G}$ a convex cone, 
by v.d.Vaart (1998) are both in terms of $\bar{\kappa}$, not $\hat\kappa$.\\
For ${\cal G}$ a (closed) convex cone, Pfanzagl+Wefelmeyer (1982) state
optimal 2-sided confidence bounds, and Janssen (1999) optimal 1-sided tests,
in terms of $\hat\kappa$, but, in the proofs, assume $- {\cal G}\subset {\cal G}$,
whence ${\cal G}$ linear, 
resp.\ $\kappa- \hat\kappa \perp {\cal G}$, whence $\hat\kappa=\bar{\kappa}$.\\ }
\end{frame}
\begin{frame}
\blu{Characterizations}
of $\bar{\kappa}\in \clin {\cal G}$ and $\hat\kappa\in \cl {\cal G}$
as in~(\ref{nnn:eq:optapprox}) by, respectively,
\begingroup \mathsurround0em\arraycolsep0em \vspace{-\smallskipamount}
\begin{eqnarray} & \Ds
  \kappa- \bar{\kappa}\perp {\cal G}\,;\hspace{.5em}\mbox{that is,\enskip}
  \langle \kappa|g \rangle = \langle \bar{\kappa}|g \rangle
  \hspace{.5em}\forall g\in {\cal G} & \\
  & \Ds
  \langle \kappa|\hat \kappa \rangle = \langle \hat{\kappa}|\hat{\kappa} \rangle
  \hspace{1em}\mbox{and}\hspace{1em}
  \langle \kappa|g \rangle \le \langle \hat{\kappa}|g \rangle
  \hspace{.5em}\forall g\in {\cal G} &
\end{eqnarray}\endgroup
\par \vspace{-\smallskipamount}Bounds based on $\hat\kappa$ are sharper since
  $   \Vert \hat \kappa \Vert < \Vert \bar{\kappa} \Vert $ unless $ \bar{\kappa}= \hat \kappa $.
\par
We shall assume either
\par \smallskip
\qquad a) ${\cal G}= \hat {\cal G}$ a closed convex cone, vertex at $0$, \qquad OR \\ 
\qquad b) ${\cal G}= \bar {\cal G}$ a closed linear space. 
\par \smallskip
For comparison, let 
       $ P=P_0\in \hat {\cal P}\subset  \bar{ {\cal P}} $, where
the smaller model $\hat {\cal P}$ has tangent set
a closed convex cone $ \hat {\cal G}$,
and the tangent set of the larger model $\bar {\cal P}$
is the closed linear span $ \bar {\cal G}= \clin \hat {\cal G}$ of~$\hat {\cal G}$.
We assume that \vspace{-\smallskipamount}
\begin{equation}
    \bar{\kappa}\in \bar{{\cal G}}\setminus \hat {\cal G}
    \hspace{.5em}\mbox{(i.e.\ $ \bar{\kappa}\ne \hat\kappa $)}\,
    \hspace{.5em} \mbox{and $\hat \kappa\ne0$\,.}
\end{equation}
\par \smallskip {\fz \blu{\bf 5.1 Example 1}\quad
Let $P= {\cal N}(0,1)$ and $\kappa(x)=x$ the identity on the real line;
$\kappa$~may be interpreted the influence curve at~$P$ of the expectation functional
as well as of the one-sample normal scores rank functional.
\par \fz
As tangent sets at~$P$, consider $\hat {\cal G}$ and $\bar{{\cal G}}$, the convex hull and
linear span, respectively, of the two tangents $g_1(x)=\zi(x)$ 
and $g_2(x)=\mu \hspace{3\ei}\zi(x)\Jc_{(|x|\le a)}$, with $a$ and $\mu=\mu_a$ in $(0,\infty)$
such that $\Vert g_2\Vert=\Vert g_1\Vert=\Vert \kappa\Vert=1 $.
\par \fz By a minimization w.r.t.~$a\in (0,\infty)$, it may be achieved that 
     $\Vert \hat \kappa \Vert = .85 \,\Vert \bar{\kappa} \Vert$.\\ }
\end{frame}
\subsection{5.2 One-Sided Tests}\begin{frame}
Given $P\in {\cal P}$, the $n$ i.i.d. observations $x_i \sim Q_n=P _{g, t/\!\sqrt{\!n}\,}$,
$n\ge1$, for any $t\in (0,\infty)$, any tangent $g\in {\cal G}$ at~$P$,
\blu{one-sided hypotheses about~$Q_n$} are
\par \medskip
\makebox[1.25em][l]{$J ^0$}:\hspace{.75em}$Q_n=P$\enskip $\iff \, g=0$ \hfill
\mbox{\sz and, employing the functional~$T$,}
\par \smallskip
\makebox[1.25em][l]{$J$}:\hspace{.75em}$\lim_{n\to \infty}
                 \sqrt{n}\,\bigl(\hspace{6\ei}T(Q_n)-T(P)\bigr) =0$\hskip.75em
                 $\iff \langle \kappa|g \rangle=0$
\par \smallskip
\makebox[1.25em][l]{$H$}:\hspace{.75em}$\lim_{n\to \infty}
                 \sqrt{n}\,\bigl(\hspace{6\ei}T(Q_n)-T(P)\bigr) \le0$\hskip.75em
                 $\iff \langle \kappa|g \rangle \le 0$
\par \smallskip
\makebox[1.25em][l]{$K$}:\hspace{.75em}$\lim_{n\to \infty}
                 \sqrt{n}\,\bigl(\hspace{6\ei}T(Q_n)-T(P)\bigr)\ge c$\hskip.75em
                 $\iff \langle \kappa|g \rangle \ge c \in (0,\infty)\:\mbox{\sz fixed}.$
\par \medskip \nz
In case $P\in \hat{\cal P}\subset \bar {\cal P}$ and corresponding tangent sets
    $\bar {\cal G}= \clin \hat {\cal G}$, the corresponding hypotheses obviously satisfy
    $ J^0 \subset \;\hat{\!\!J} \subset \;\bar{\!\!J}$, $\,\hat{\!H}\subset \,\bar{\!H}$,
    $\,\hat{\!K}\subset \,\bar{\!K} $.
\par  Depending on the choice ${\cal G}=\hat{{\cal G}}\; \mbox{or}\;\bar{{\cal G}}$, put
    $\tilde \kappa= \hat{\kappa}$, respectively $\tilde \kappa=\bar{\kappa}$.
\par We consider sequences $\tau=(\tau_n)$ of tests $\tau_n$ at sample size~$n$.
\par \medskip
\blu{\bf 5.2 Theorem 1}\enskip \blu{[\,$J^0$ vs.~$K$\,]}\quad {\sl If \enskip
  $ \limsup _{n\to \infty} \int \tau_n \,dP^n \le \alpha $ \enskip then
\begin{equation}
    \inf _{K}\limsup _{n\to \infty} \smash{\int }\tau_n \,dQ_n^n \le 
    \Phi \big(- u _{\alpha}+ \Tfrac{c}{\Vert \tilde \kappa \Vert}\,\bigr)
\vspace{- \smallskipamount}\end{equation}
The power bound is achieved uniquely---up to $\Lo _{P^n}(n^0)$---by the tests
\begin{equation}
    \tilde {\tau}_n = \lJc \bigl\{ {\Ts
    \frac{1}{\sqrt{n}\,}\sum _{i=1}^n }
    \tilde \kappa (x_i) > \Vert \tilde \kappa \Vert u _{\alpha}\bigr\}
\vspace{-\smallskipamount}\end{equation}\/}
\blu{[\,$\,\bar{\!H}$ vs.~$\bar{K}$\,]}\quad {\sl
In case ${\cal G}= \bar{{\cal G}}$ moreover\quad
  $ \sup _{\,\bar{\!H}}\limsup _{n\to \infty} \int \bar\tau_n \,dQ_n^n \le \alpha $\/}
\end{frame}
\begin{frame}\par \smallskip {\sz
\blu{\bf Proof}\enskip The closed convex set $G_1=$ all $g_1\in {\cal G}$ such that
   $\langle \kappa|g_1 \rangle \ge c$ has minimum norm element
   $q_1= \tilde t \tilde \kappa$ with $\tilde t = c/\Vert \tilde \kappa \Vert^2$.
Thus 3.4 Theorem~1 provides the unique asymptotic maxmin test $\tilde \tau$ for $J^0$ vs.~$K$.
To enlarge the null $J^0$ to $J$ or $H$, set $G_0=$ all $g_0\in {\cal G}$ such that
   $\langle \kappa|g_0 \rangle = 0$, respectively $\le 0$.
\par
In case ${\cal G}= \bar{{\cal G}}$, $q_1=q _{10}=q_1-q_0$, with $q_0=0$, turns out of
minimum norm also in $\,\bar{\!G}_{10}= \,\bar{\!G}_1- \,\bar{\!G}_0$. This is true since
  $ c \le \langle \bar{\kappa}|g_1 \rangle - \langle \bar{\kappa}|g_0 \rangle \so
    \Vert q _{10} \Vert^2 \le \langle q _{10}|g _{10} \rangle$ for all
  $g _{10}\in \,\bar{\!G}_{10}$, and thus~(\ref{nnn:eq:optapprox}).
  3.4 Theorem~1 now applies again for $\,\bar{\!H}$ vs.~$\,\bar{\!K}$.
\par \sz
In case ${\cal G}= \hat{{\cal G}}$, minimization of the norm 
in $\,\hat{\!G} _{10}=\,\hat{\!G} _{1}-\,\hat{\!G} _{0}$ is \pur{ yet unsolved}. }
\par \medskip
For $\bar{{\cal G}}=\clin \hat{{\cal G}}$ note that $\hat{\kappa}\ne \bar{\kappa}$
iff $\langle \kappa|g \rangle  < \langle \hat \kappa|g \rangle $ for some
    $g\in \hat{{\cal G}}$.
\par \smallskip
\blu{\bf 5.2 Theorem 2}\enskip {\sl In case ${\cal G}= \hat{{\cal G}}$ 
assume some $g_0\in \hat{{\cal G}}$ such that \vspace{-\medskipamount}
\begin{equation}\label{ccc:eq:cond<}
   \langle \kappa|g_0 \rangle \le 0 < \langle \hat\kappa|g_0 \rangle
\vspace{-\smallskipamount}\end{equation}
Then\hskip.75em 
$ \sup _{\;\hat{\!\!J}} \liminf _{n\to \infty} \int \hat{\tau}_n\,dQ_n = 1$}\enskip  
\pur{\sf [\,level breakdown of $\hat{\tau}$ on $\;\hat{\!\!J}$\:]}
\par \smallskip {\fz
\blu{\bf 5.2 Example 3}\quad In 5.1 Example 1, although $\hat{\kappa}\ne \bar{\kappa}$,
condition~(\ref{ccc:eq:cond<}) is not fulfilled.\\
But replace tangent $g_2$ there by
  $g_3(x)=\delta \Jc _{(0,a]}(x)- \eta \Jc _{(a,\infty)}(x)=-g_3(-x) $, $x\ge0$,
where the constants may be determined such that $\Vert g_3 \Vert=1$.
Then~$g_3$ achieves~(\ref{ccc:eq:cond<}). }
\par \smallskip {\sz
\blu{ {\bf 5.2 Remark 4} [\,$\bar{\tau}$ for $\,\hat{\!H}$ vs.~$\hat{\!K}$\,]}\enskip {\sl 
In case $\bar{{\cal G}}=\clin \hat{{\cal G}}$, since $P\in \,\hat{\!H}\subset \,\bar{\!H}$,
  the test $\bar{\tau}$ achieves
  $ \sup _{\,\hat{\!H}} \mbox{asy.level of $\bar{\tau}_n$}
    = \alpha  $
and, with $\inf _{\,\hat{\!K}}$ attained 
at $ \hat{q}_1 = \hat{t}\hat{\kappa}\in \hat{\!K} \subset \,\bar{\!K}$,
\begin{equation}\Ts
    \inf _{\,\hat{\!K}} \mbox{asy.power of $\bar{\tau}_n$}
    = \Phi \big(- u _{\alpha}+ \Tfrac{c}{\Vert \bar\kappa \Vert}\,\bigr) \: \mbox{\fz (\,$\bar{\tau}$ best?)}\:
    < \Phi \big(- u _{\alpha}+ \Tfrac{c}{\Vert \hat\kappa \Vert}\,\bigr)
\end{equation}\null\/}}
\end{frame} 
\subsection{5.3 One- and Two-Sided Confidence Bounds}\begin{frame} 
Given $P$, ${\cal G}$,  
$T$ differentiable under 
  $P _{n,g,t}:=P _{g, t/\!\sqrt{n}\,}$,  
as in (\ref{ccc:eq:Pgs}), (\ref{ccc:eq:Tkappa}). 
\par  
Consider estimator sequences $S=(S_n)$ which, for certain 
tangents~$g\in {\cal G}$, asymptotically 
have median $\ge T$ or $\le T$ such that, respectively, $\forall\,t>0$, 
\vspace{-\smallskipamount}\begingroup \mathsurround0em\arraycolsep0em
\begin{eqnarray}\label{ccc:eq:asymed+g} 
   \limsup _{n\to \infty} P^n _{n,g,t} \bigl\{ 
   S_n < T(P _{n,g,t})    \bigr\} &{}\le {}& \Tfrac{1}{2} 
\\ \noalign{\noindent \vskip-\medskipamount\nopagebreak} 
\label{ccc:eq:asymed-g}
   \limsup _{n\to \infty} P^n _{n,g,t} \bigl\{
   S_n>T(P _{n,g,t})   \bigr\}    &{}\le {}& \Tfrac{1}{2} 
\end{eqnarray}\endgroup
\par \vspace{-1.5 \smallskipamount}
We assume ${\cal G}$ closed, a) a convex cone $\hat{{\cal G}}$, or
b)  a linear space $\bar{{\cal G}}$.
\par \smallskip 
\blu{ {\bf 5.3 Theorem 1}\quad a) ${\cal G}=\hat{{\cal G}}$:}\enskip 
{\sl If\/ {\sf (\ref{ccc:eq:asymed+g})} holds for $g= \hat{\kappa}$, then 
  $\:\forall\, c>0$ 
\vspace{-.75\smallskipamount}\begin{equation}
   \limsup _{n\to \infty} P^n \bigl\{ T(P)>S_n- \Tfrac{c}{\sqrt{n}\,}\bigr\} 
   \le \Phi \bigl( \Tfrac{c}{\Vert \hat{\kappa} \Vert} \bigr) 
\vspace{-1.5\smallskipamount}\end{equation} 
The upper bound is attained by $\,\hat{\!S}\hskip.75em \forall\, c>0$, 
iff 
\vspace{-\smallskipamount}\begin{equation}\label{ccc:eq:hutS}\Ts
   \sqrt{n}\,\bigl(\,\hat{\!S}_n-T(P)\bigr)_+ = 
   \bigl( n ^{-1/2}\Tsum _{i=1}^n \hat{\kappa}(x_i)\bigr)_+ + \Lo _{P^n}(n^0)
\end{equation}
\par \blu{\sf b) ${\cal G}=\bar{{\cal G}}$:}\enskip 
Under\/ {\sf (\ref{ccc:eq:asymed+g})} for $g= \bar{\kappa}$ and\/ 
             {\sf (\ref{ccc:eq:asymed-g})} for $g= -\bar{\kappa}$, 
   then $\;\forall\,c',c''>0$, 
\vspace{-.75\smallskipamount}\begin{equation}
   \limsup _{n\to \infty} P^n \bigl\{ 
   S_n-\Tfrac{c''}{\sqrt{n}\,}<T(P)<S_n+\Tfrac{c'}{\sqrt{n}\,}\bigr\}
   \le \Phi \bigl( \Tfrac{c''}{\Vert \bar{\kappa} \Vert}\bigr) - 
   \Phi \bigl( \Tfrac{-c'}{\Vert \bar{\kappa} \Vert}\bigr)
\vspace{-1.5\smallskipamount}\end{equation}  
The upper bound is attained by $\,\bar{\!S}\hskip.75em \forall\,c',c''>0$ 
iff
\vspace{-\smallskipamount}\begin{equation}\label{ccc:eq:barS}\Ts
   \sqrt{n}\,\bigl(\,\hat{\!S}_n-T(P)\bigr) =
   n ^{-1/2}\Tsum _{i=1}^n \bar{\kappa}(x_i) + \Lo _{P^n}(n^0)
\end{equation}\null\/}
\end{frame}
\begin{frame} Estimators such that, with any $\eta\in L_2(P)$, $\eta \perp 1$,  
\vspace{-\smallskipamount}\begin{equation} 
  \sqrt{n}\,\bigl(S_n-T(P)\bigr) =
   n ^{-1/2}\Tsum _{i=1}^n \eta(x_i) + \Lo _{P^n}(n^0)
\vspace{-\smallskipamount}\end{equation}
for all tangents~$g\in {\cal G}$, all $t>0$, are asymptotically normal 
\vspace{-\smallskipamount}\begin{equation}
  \Lw _{P^n _{n,g,t}} \bigl\{ \sqrt{n}\,\bigl(S_n-T(P_{n,g,t})\bigr) \bigr\} 
  \longrightarrow {\cal N}\bigl( t \langle \eta-\kappa|g \rangle, \Vert \eta \Vert^2\bigr)
\end{equation} 
\par \blu{{\bf 5.4 Corollary 2}\quad 
            [\,${\cal G}=\bar{{\cal G}}$, stability of $\,\bar{\!S}$\,]:}\enskip  
{\sl The estimator $\bar{S}$ achieves 
\vspace{-\smallskipamount}\begin{equation} 
   P^n _{n,g,t} \bigl\{
   \,\bar{\!S}_n-\Tfrac{c''}{\sqrt{n}\,} <T(P_{n,g,t}) 
   < \,\bar{\!S}_n+\Tfrac{c'}{\sqrt{n}\,}    \bigr\}
   \longrightarrow \Phi \bigl( \Tfrac{c''}{\Vert \bar{\kappa} \Vert} \bigr) 
   -  \Phi \bigl( \Tfrac{-c'}{\Vert \bar{\kappa} \Vert} \bigr)
\vspace{-\smallskipamount}\end{equation}
for all $g\in \bar{{\cal G}}$, $t>0$, and all $c',c''\ge0$; in particular, 
is asy.\,median unbiased achieving $\lim_n=\frac{1}{2}$ 
in {\sf (\ref{ccc:eq:asymed+g}), (\ref{ccc:eq:asymed-g})} 
  $\forall\, g\in \bar{{\cal G}}$.\/} 
\par \smallskip 
In case ${\cal G}= \hat{{\cal G}}\subset \bar{{\cal G}}$ and $\hat{\kappa}\ne \bar{\kappa}$, 
   $\exists\, g_1\in \hat{{\cal G}}$ such that 
   $0<\langle \kappa|g_1 \rangle< \langle \hat{\kappa}|g_1  \rangle$. 
Consequently, 
no opt.~estimator $\,\hat{\!S}$ of form~{\sf (\ref{ccc:eq:hutS})} may fulfill 
condition~{\sf (\ref{ccc:eq:asymed-g})}. Moreover,  as lower confidence limit, 
   $\,\hat{\!S}$ \pur{breaks down} under~$P_{n,g_1,t}$. 
\par \smallskip 
\blu{{\bf 5.4 Proposition 3}\quad
     [\,${\cal G}=\hat{{\cal G}}$, positive asy.\,bias of $\,\hat{\!S}$\,]:}\enskip 
{\sl  $\exists\,g_1\in \hat{{\cal G}}$ such that any optimal 
estimator $\,\hat{\!S}$ of form~{\sf (\ref{ccc:eq:hutS})} satisfies, $\forall\,c\ge0$,  
\vspace{-.75\smallskipamount}\begin{equation}
   \lim _{t\to \infty}\lim _{n\to \infty}
   P^n_{n,g_1,t} \bigl\{ T(P_{n,g_1,t})\ge S_n- \Tfrac{c}{\sqrt{n}\,}\bigr\}
   =0 < \Phi \bigl( \Tfrac{c}{\Vert \hat{\kappa} \Vert} \bigr) 
\vspace{-.75 \smallskipamount}\end{equation}
in particular, violates~{\sf (\ref{ccc:eq:asymed-g})} as 
  $  \lim _{t}\lim _{n} 
     P^n _{n,g_1,t} \bigl\{\,\hat{\!S}_n>T(P _{n,g_1,t})   \bigr\} =1 >\Tfrac{1}{2}
  $.\/}  
\end{frame}
\section{6. Unknown Neighborhood Radius}
\subsection{6.1 List of Ideal Models}\begin{frame}\blu{\bf Models:}\hfill
\mbox{\tiny Rieder, Kohl, Ruckdeschel (2008)}
\begin{itemize}
\item Location: $y = \theta + u$, $u\sim{\cal N}_k(0, \EM_k)$,
      $P_0 = {\cal N}_k(0, \EM_k)=P$
\item Scale ($k=1$): $y = \sigma u$, $u\sim{\cal N}(0,1)=P_1=P$
\item Regression ($k \ge 1$): \qquad $y = x\,\theta + u$;\quad
      $x$, $u$ sto.~indep.
\vspace{-\smallskipamount}\begin{eqnarray*}
  &u\sim {\cal N}(0,1)\;,\; x\sim K(dx)\\
   & P=P_0(dx, du) = K(dx)\,{\cal N}(0,1)(du)
\end{eqnarray*}
\par \vspace{-\medskipamount}%
For $\alpha=2$, coincidence of results with 1-dim.~location.
\par For $\alpha=1$, 
assume $K$ spherically (elliptically) symmetric.
\item   ARMA($p,q$)-models (with shift) are covered, setting~$K=\Lw(H)$.\\
Ideal innovations i.i.d.~$\sim {\cal N}(0,1)$, then $K$ multivariate normal.\\
For $\alpha=2$, coincidence of results with 1-dim.~location.
\item ARCH(1): $y_t= \sqrt{1+\theta y _{t-1}^2}\:u_t$,
      $u_t$ i.i.d.~$\sim {\cal N}(0,1)$\\ 
For $\alpha=2$, coincidence of results with 1-dim.~scale.
\end{itemize}
\end{frame}
%
%
\subsection{6.2 Robust Neigborhoods}\begin{frame}
\blu{\bf Neigborhoods About \boldmath $P$ \unboldmath:}
\medskip
\begin{itemize}
\item (1-dim.~location) 
symmetric contamination~nbd of size $s\in [0,1)$ :
\begin{center}
  $ F = (1 - s)\,{\cal N}(0,1) + s\,H,\quad H \mbox{\ symmetric}$
\end{center}
\item   $r/\!\sqrt{n}$ - nbds at sample size $n$:
\begin{center}
  $Q_n = (1 - \frac{r}{\sqrt{n}})P + \frac{r}{\sqrt{n}}\,H$
\end{center}
(location, unconditional regression; \quad scale: H symmetric)
\item conditional regression $r/\!\sqrt{n}$\,-\,nbds,
      with radius curve $\varepsilon (x)$: 
\begin{center}
  $Q_n(du\,|\,x) = (1 - \frac{r}{\sqrt{n}}\,\varepsilon(x))\Phi(du)
  + \frac{r}{\sqrt{n}}\,\varepsilon(x)\,H(du\,|\,x)$,
\end{center}
\item in time series: contaminated transition probabilities\\
   $\Ds Q_n(dy_t|\bar y_{t-1})$ \hspace{5em}
   where $\bar y_{t-1} : = y_{t-1},\ldots, y_1$
   $$\Ts {} = (1 - \frac{r}{\sqrt{n}}\,\varepsilon(\overline{y}_{t-1}))
       P(dy_t|\bar y_{t-1}) + \frac{r}{\sqrt{n}}\,\varepsilon(\overline{y}_{t-1})\,
       H_n(dy_t\;|\;\overline{y}_{t-1})$$
\item $\|\varepsilon \|_\alpha\le 1$:\quad $\Ew \varepsilon\le 1$
($\alpha=1$),\quad $\Ew \varepsilon^2\le 1$
($\alpha=2$),\quad $\varepsilon\le 1$ ($\alpha=\infty$)\\
  $\Ew$ is taken under the ideal measure $P$, resp.\ ideal regressor distr.
\end{itemize}
\end{frame}
%
%
\subsection{6.3 Relative Maximum Risk}\begin{frame}
\blu{\bf Relative Maximum Risk Over Neighborhoods:}
\par \medskip
We use the estimate which is optimally robust for the neighborhood
model of an assumed radius while this radius may not be true.
\begin{itemize}
\item relative Var (in Huber[64] model):
\par (minmax) M-estimates of location,
      $ \sum \limits_{i=1}^n \psi(y_i-S_n) \approx 0$
$${\rm relVar}\,(\psi_{s_0}, s) =
   \frac{{\rm maxVar}\,(\psi_{s_0},s)} {{\rm maxVar}\,(\psi_s, s)},
    \qquad 0\le s < 1$$
\item relative MSE ($r/\!\sqrt{n}$ - neighborhoods, Ri[94]):
\par (minmax) asy.~linear estimates with influence curves
\begin{center}
  $\Ts \sqrt{n}\,(S_n- \theta) - n ^{-1/2} \sum
    \limits_{i=1}^n \eta(y_i) \longrightarrow 0 $\quad in~$P$-prob.
\end{center}
$${\rm relMSE}\,(\eta_{r_0}, r) =
    \frac{{\rm maxMSE}\,(\eta_{r_0},r)}{{\rm maxMSE}\,(\eta_r, r)},
    \qquad 0\le r<\infty$$
\end{itemize}
\end{frame}
%
%
\subsection{6.4 Location (1-dimensional)}\begin{frame}
\blu{\large\bf Location (1-dim)}
\par \bigskip
\blu{Minimax asymptotic variance} 
\medskip
\begin{itemize}
\item Minimax M-estimate for $s\in [0,1)$:
\smallskip\begin{center}
  $\psi_{s}(u) = (-m_s)\lor  u\wedge m_s$,\qquad
  $\frac{s}{1-s}\;m_s = \Ew (|u| - m_s)_+$
\end{center}\smallskip
\item Maximal asymptotic variance of $\psi_{s_0}$ under $s$:
\smallskip\begin{center}
  $\Ds {\rm maxVar}\,(\psi_{s_0}, s) =
  \frac{(1-s)\Ew \psi^2_{s_0} + s\,m^2_0}
  {\bigl[(1-s)\Ew \psi'_{m_0}\bigr]^2}\;$
\end{center}\medskip
\item Median ($s = 1$): $\psi_1(u) = \sign{(u)} =
\lim\limits_{s\to 1}\frac{1}{m_s}\psi_s(u)$,
$$\Ts {\rm maxVar}\,(\psi_1, s) = \frac{\pi}{2(1-s)^2},
  \qquad{\rm relVar}\,(\psi_1, s)\longrightarrow 1\quad (s\to 1)$$
\end{itemize}
\end{frame}
%
%
\begin{frame}
\par \smallskip
\blu{Minimax asymptotic MSE} 
\medskip
\begin{itemize}
\item Minimax IC for $r \in [0, \infty)$:\quad
  $\eta_{r}(u) = A_{r}u \min{ \left\{ 1,\, \frac{c_{r}}{|u|}\right\}}$, 
\begin{center}
  $1 = A_{r}\Ew u^2
  \min{\left\{1,\,\frac{c_{r}}{|u|}\right\}}$,
  $\qquad r^2\,c_{r} = \Ew \left(|u| - c_{r}\right)_+$
\end{center}
\smallskip\item Median ($r = \infty$):
  $\eta_\infty\,(u) = b_{{\rm min}}\sign{(u)}$,
\smallskip\item Minimal bias (of ALE):
  $b_{{\rm min}}=\sqrt{\frac{\pi}{2}}$
\smallskip\item Maximal MSE of $\eta_{r_0}$ under $r$:
\smallskip\begin{center}
  ${\rm maxMSE}\,(\eta_{r_0}, r) = A^2_{r_0}\;\Ew
  \min{\left\{u^2,\,c^2_{r_0}\right\}}
  \; +\; r^2\;A^2_{r_0}c^2_{r_0}$
\end{center}\end{itemize}
\par \bigskip \pur{Coincidence}\vspace{-\smallskipamount}
\begin{eqnarray*}\quad
  (1-s)\,{\rm maxMSE}\,(\eta_{r_0}, r) &=& {\rm maxVar}\,(\psi_{s_0}, s)\\
  \quad\so   {\rm relVar}\,(\psi_{s_0}, s) &=& {\rm relMSE}\,(\eta_{r_0}, r)
\end{eqnarray*}
where $r$ and $s$ correspond via \enskip $s = r^2/(1+r^2)$.
\end{frame}
\begin{frame}\vspace{-5cm}\begin{center}
\includegraphics[width=0.8\linewidth]{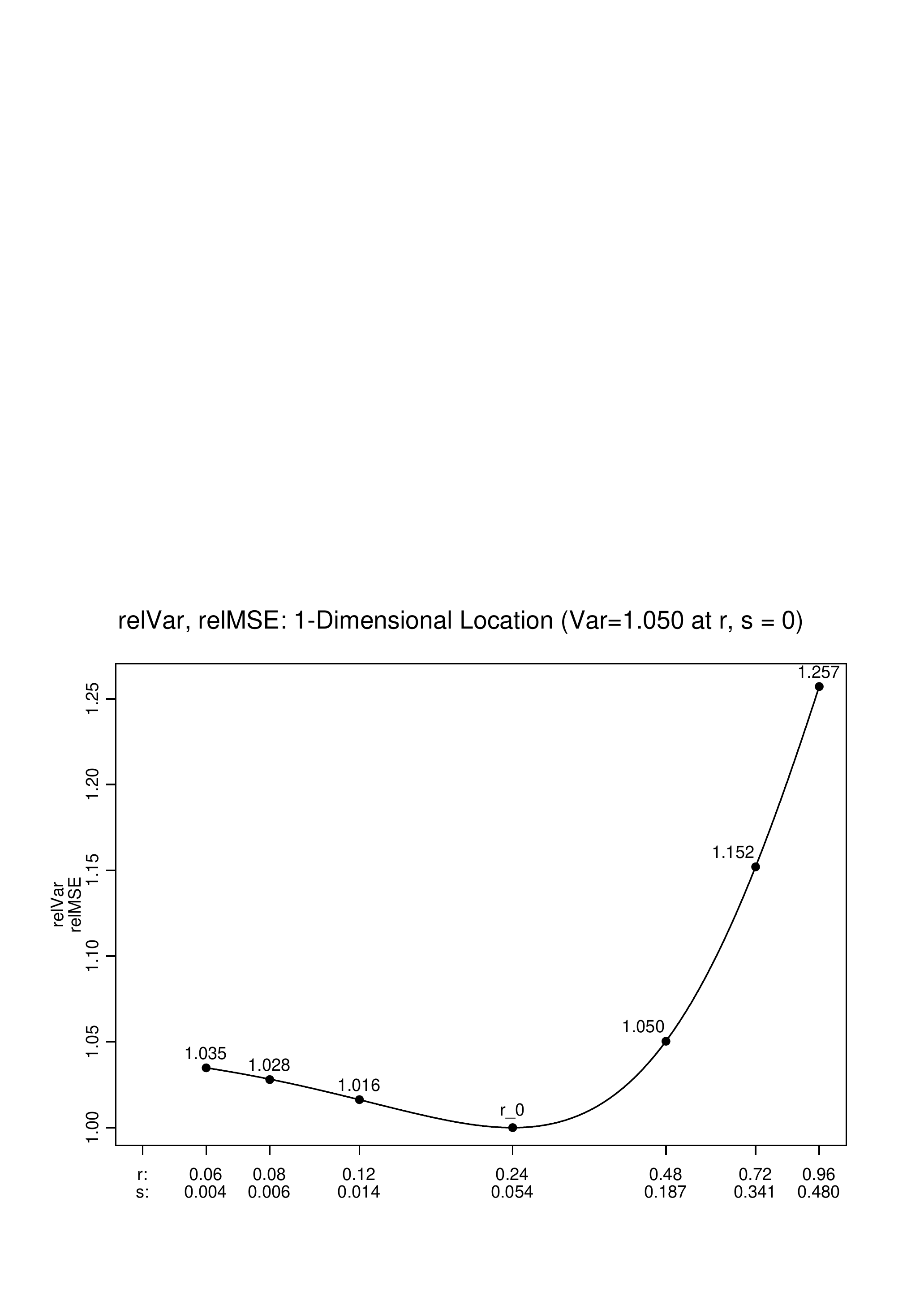}
\end{center}
\end{frame}
\begin{frame}\vspace{-5cm}\begin{center}
\includegraphics[width=0.8\linewidth]{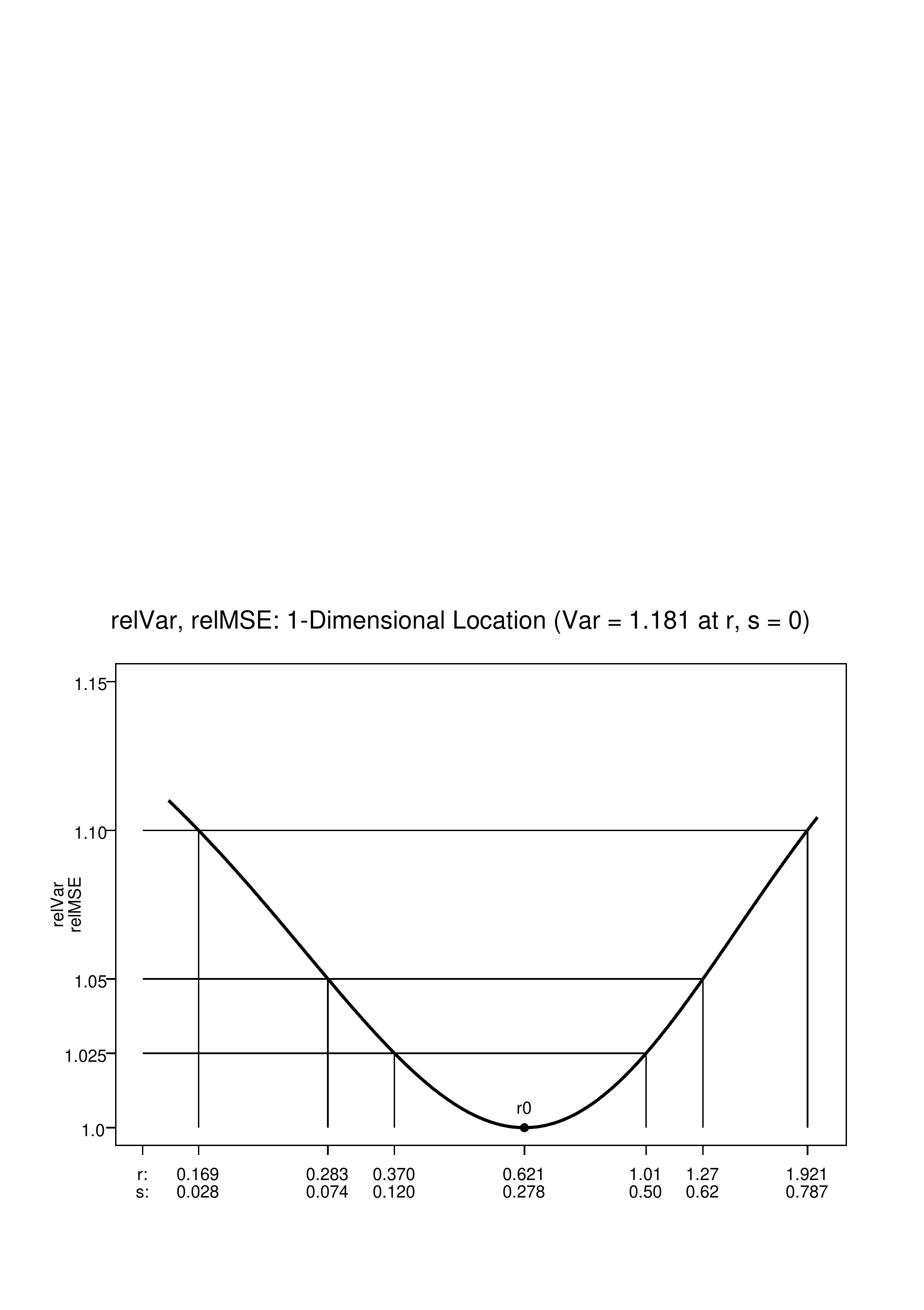}
\end{center}
\end{frame}
\begin{frame}\vspace{-5cm}\begin{center}
\includegraphics[width=0.8\linewidth]{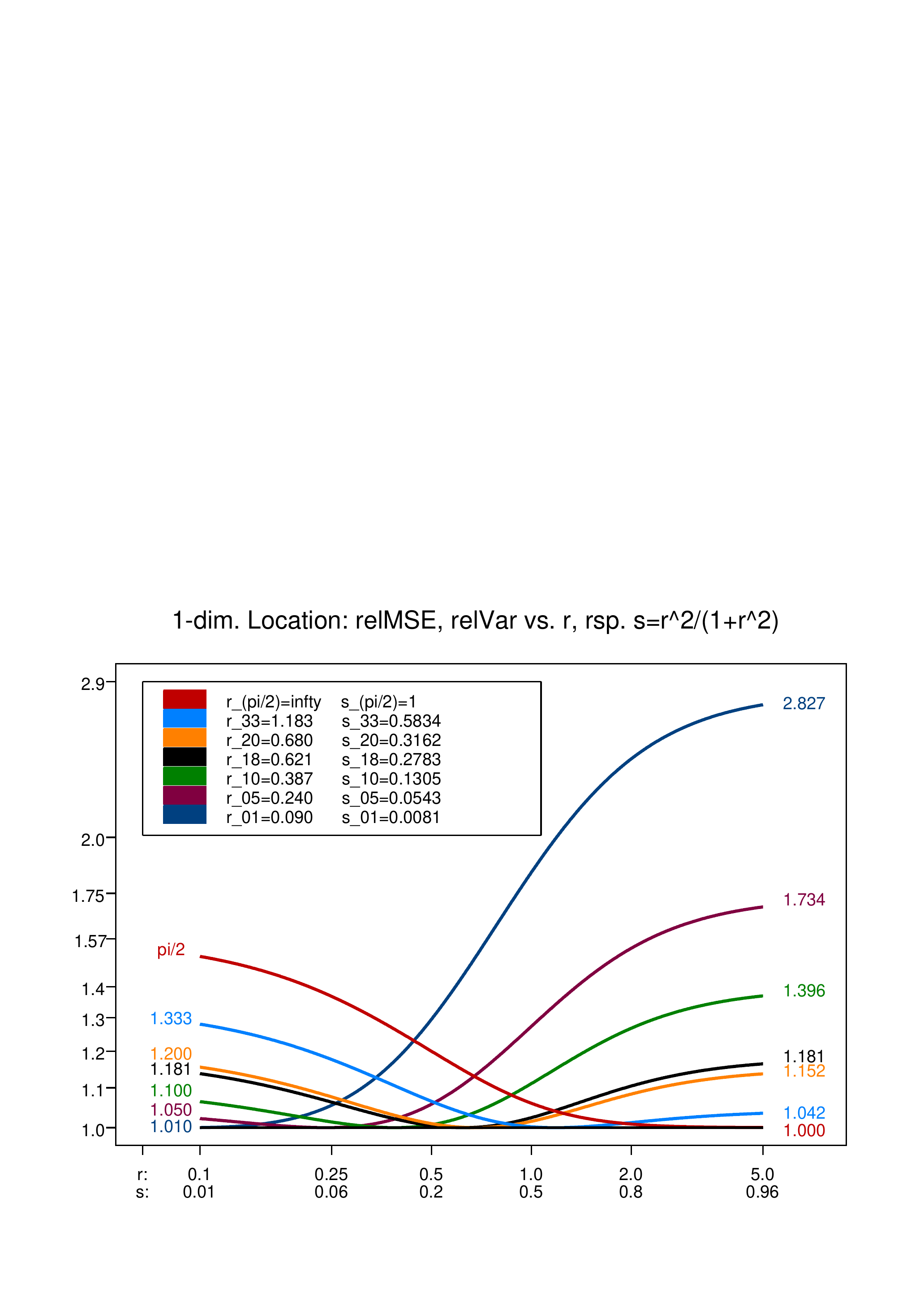}
\end{center}
\end{frame}
%
%
\subsection{6.5 Location ($k$-dimensional)}\begin{frame}
\blu{\large\bf Location (k-dim)}
\par \medskip \blu{Minimax asymptotic MSE}
\medskip
\begin{itemize}
\item Minimax IC for $r \in [0, \infty)$:\quad
  $\eta_{r}(u) = \alpha_{r}u   \min{\left\{1, \frac{c_{r}}{|u|}\right\}}$, 
\begin{center}
  $k = \alpha_{r}\,\Ew |u|^2
  \min{\left\{1,\,\frac{c_{r}}{|u|}\right\}}$,
  $\quad r^2\, c_{r} = \Ew (|u| - c_{r})_+$
\end{center}
\smallskip\item min-$L_1$ ($r = \infty$):\hspace{.5em}
  $\sum_{i=1}^n|u_i - \hat\theta|_2 = \min_\theta{}!$\hspace{.75em}
  $\eta_\infty\,(u) = b_{{\rm min}}\frac{u}{|u|}$
\smallskip
\item  Minimal bias (of ALE):
\begin{center}
  $\Ds b_{{\rm min}} = \frac{k}{\Ew |\Lambda|} =
  \frac{k\,\Gamma(\frac{k}{2})}{\sqrt{2}\,\Gamma(\frac{k+1}{2})}$,\quad
  ${\Ds\quad\frac{b_{{\rm min}}}{\sqrt{k}}\to 1}$,
  ${\Ds\frac{\Ew |\eta_\infty|^2}{k}\to 1}$
\end{center}
\smallskip\item Maximal MSE of $\eta_{r_0}$ under $r$:
\smallskip\begin{center}
  ${\rm maxMSE}\,(\eta_{r_0}, r) = \alpha^2_{r_0}\Ew
  \min{\{|u|^2,\,c^2_{r_0}\}}
  \; +\; r^2\;\alpha^2_{r_0}c^2_{r_0}$
\end{center}
\bigskip\item Relative MSE: \quad \pur{$\eta _{\infty}$~becomes radius--minimax}
$$ \Nlim _{k\to \infty}\:
  \frac{{\rm maxMSE}\,(\eta_{r_0}, r)}{{\rm maxMSE}\,(\eta_\infty, r)}
  = 1 
$$
uniformly in $0 \le r_0,r \le \mbox{any }r_1 < \infty$.
\end{itemize}
\end{frame}
\begin{frame}\vspace{-5cm}\begin{center}
\includegraphics[width=0.8\linewidth]{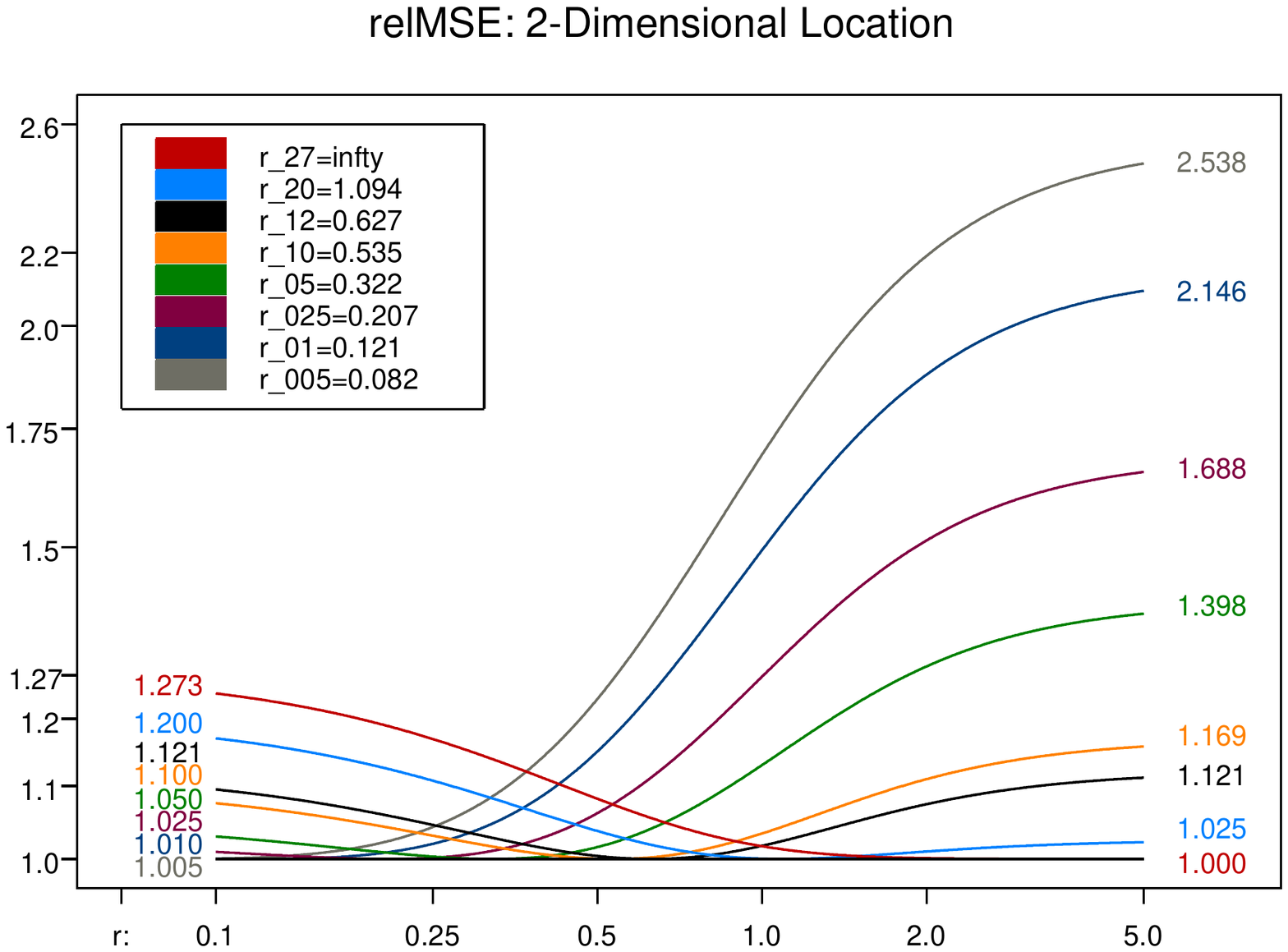}
\end{center}
\end{frame}
%
%
\subsection{6.6 Regression ($k$-dimensional)}\begin{frame}
\blu{\large\bf Regression (k-dim)}
\par \medskip \blu{Minimax asymptotic MSE $(*=c, \alpha=1)$}
\begin{itemize}
\item Minimax IC for $r \in [0, \infty)$:
  $\eta_{r}(x, u) = \alpha_{r}x\,u
  \min{\left\{1, \frac{c_{r}}{|x u|}\right\}}$, 
\begin{center}
  $k = \alpha_{r}\,{\rm E}\,|x|^2 u^2
  \min{\left\{1,\,\frac{c_{r}}{|x u|}\right\}}$,
  $\quad r^2\, c_{r} = {\rm E}\,(|x u| - c_{r})_+$
\end{center}
\smallskip\item weighted min-$L_1$ ($r=\infty$):\quad
$\eta_\infty\,(x, u) = b_{{\rm min}}\frac{x}{|x|}\sign{(u)}$
\item Minimal bias (of ALE):\quad
  $b_{{\rm min}} = \frac{k}{{\rm E}\,|\Lambda|}
  =\sqrt{\frac{\pi}{2}}\,\frac{k}{{\rm E}\,|x|}$
\smallskip\item Maximal MSE of $\eta_{r_0}$ under $r$:
\smallskip\begin{center}
  ${\rm maxMSE}\,(\eta_{r_0}, r) =
  \alpha^2_{r_0}\Ew \min{\{|x|^2u^2,\,c^2_{r_0}\}}
  \; +\; r^2\;\alpha^2_{r_0}c^2_{r_0}$
\end{center}
\medskip\item RelMSE same for all~$\theta$, but depends on~$K(dx)$.
\par \medskip \pur{Convergence to 1-dim. location}
$$ \Nlim _{k\to \infty}{\rm relMSE}\,(\eta_{r_0}, r)=
  {\rm relMSE}\,(\eta^{\rm 1loc}_{r_0},  r)
$$
uniformly for $0 \le r_0,r \le \mbox{any }r_1 < \infty$, as $k\to \infty$.\\
In case $(*=c, \alpha=2)$ limit attained $\forall\,k\ge1$.
\end{itemize}
\end{frame}
\begin{frame}\vspace{-5cm}\begin{center}
\includegraphics[width=0.8\linewidth]{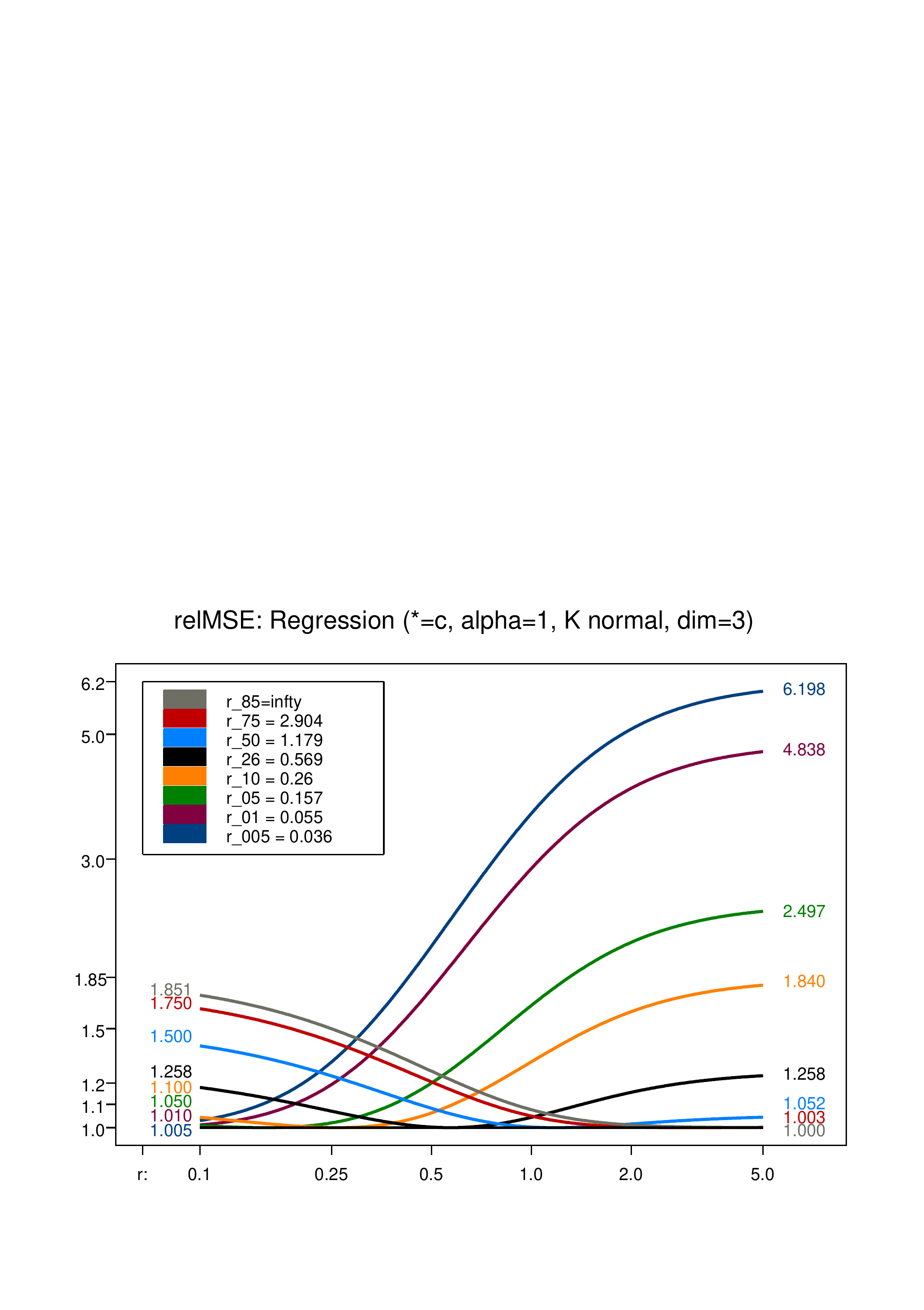}
\end{center}
\end{frame}
%
%
\subsection{6.7 Scale ($1$-dimensional)}\begin{frame}
\blu{\large\bf Scale (1-dim)}
\par \medskip \blu{Minimax asymptotic MSE 
                        for $(*=c)$ contamination balls}\sz
\begin{itemize}
\item Minimax IC for $r \in [0, \infty)$:\hspace{1.5em}%
  $\Ds \eta_{r}(u) = A_{r}(u^2-\alpha_r^2)
  \min{\left\{1, \frac{c_{r}}{|u^2-\alpha_r^2|}\right\}}$ 
\begin{eqnarray*}
  0 &=& \Ts \Ew(u^2-\alpha_r^2)
  \min{\left\{1,\,\frac{c_{r}}{|u^2-\alpha_r^2|}\right\}},\\
  A_r^{-1} &=& \Ts \Ew(u^2-\alpha_r^2)^2
  \min{\left\{1,\,\frac{c_{r}}{|u^2-\alpha_r^2|}\right\}},\\
  \qquad r^2 c_{r} &=& \Ew(|u^2-\alpha_r^2| -
  c_{r})_+
\end{eqnarray*}
\item MAD ($r = \infty$):\quad
  $\eta_\infty\,(u) = b_{{\rm min}}\sign{(|u| - \alpha_\infty)}$,\hskip.75em 
  $\hat\theta = \alpha_\infty ^{-1}\:{\med}(|u_i|)$
\smallskip\item Minimal bias (of ALE):\quad
  $b_{{\rm min}}=(4\alpha_\infty\varphi(\alpha_\infty))^{-1}= 1.166$
   \item $0<\alpha_r$ decreasing from $\alpha_0 = 1$ to
   $\alpha_\infty := \Phi^{-1}(3/4) = 0.674$
\smallskip\item \pur{ clipping of $|u|$ only from above for $r\le 0.92$;\\
      clipping of $|u|$ from below and above iff $r\ge 0.92$ }
\smallskip\item
For $r_0,r\in[0,\infty)$, the maximal MSE is
$$ 
  {\rm maxMSE}\,(\eta_{r_0}, r) = A^2_{r_0} \Ew 
  \min{\{|u^2 - \alpha_{r_0}^2|^2,\,c^2_{r_0}\}}
  \; +\; r^2A^2_{r_0}c^2_{r_0}
$$
\end{itemize}\nz
\end{frame}
\begin{frame}\blu{Minimax asymptotic MSE
                        for $(*=v)$ contamination balls}\sz 
\par \vspace{\bigskipamount}
\begin{itemize}
\item Minimax IC for $r \in [0, \infty)$:\smallskip
\begingroup \mathsurround0em\arraycolsep0em
\begin{eqnarray*} 
  \eta_{r}(u) &{}={}&
  A_{r}\{[\hspace{9\ei}g_r \lor  u^2 \wedge
  (\hspace{6\ei}g_r + c_r)] - 1\} \\
\noalign{\vspace{\medskipamount}}
  0 &{}={}& \Ew(\hspace{6\ei}g_r-u^2)_+ - \Ew(u^2-g_r-c_r)_+\\
  1 &{}={}& A_r\Ew u^2\big\{[\hspace{9\ei} 
            g_r \lor  u^2 \wedge (\hspace{6\ei}g_r + c_r)]
            - 1\big\}\\
  \qquad r^2\, c_{r} &{}={}& \Ew(\hspace{6\ei}g_r-u^2)_+
\end{eqnarray*}\endgroup
\item MADv ($r = \infty$): 
$$ 
  \eta_\infty\,(u) = \omega_v^{{\rm  min}}\hspace{6\ei}
  \big\{P(|u| < 1)\Jc(|u| > 1) - P(|u| > 1) \Jc(|u| <  1)\big\}
$$ 
\item Minimal bias (of ALE):
$$ 
  \omega_v^{{\rm min}} = (\hspace{6\ei}\Ew\Lambda_+)^{-1}
  = \sqrt{\Tfrac{\pi}{2}\,e}\approx 2.066
$$ 
\item \pur{clipping of $|u|$ always from above and below}
\item
For $r_0,r\in[0,\infty)$, the maximal MSE is
$$ 
  {\rm maxMSE}\,(\eta_{r_0}, r) = A^2_{r_0}\Ew 
  \big\{[\hspace{9\ei}g_{r_0} \lor  u^2 \wedge (\hspace{6\ei}
  g_{r_0} + c_{r_0})] - 1\big\}^2 + r^2A^2_{r_0}c^2_{r_0}
$$ 
\end{itemize}\nz
\end{frame}
\begin{frame}\vspace{-5cm}\begin{center}
\includegraphics[width=0.8\linewidth]{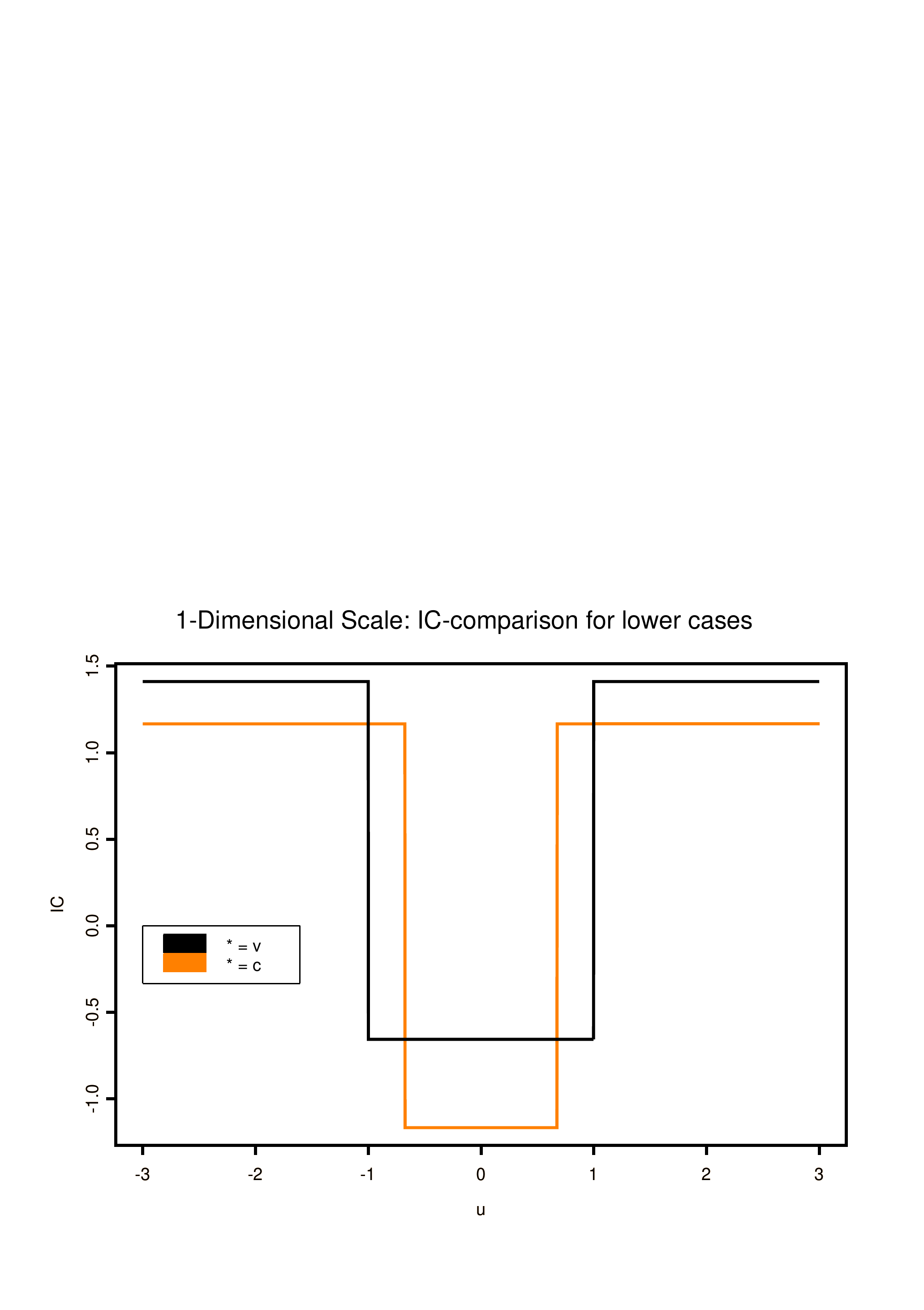}
\end{center}
\end{frame}
\begin{frame}\vspace{-5cm}\begin{center}
\includegraphics[width=0.8\linewidth]{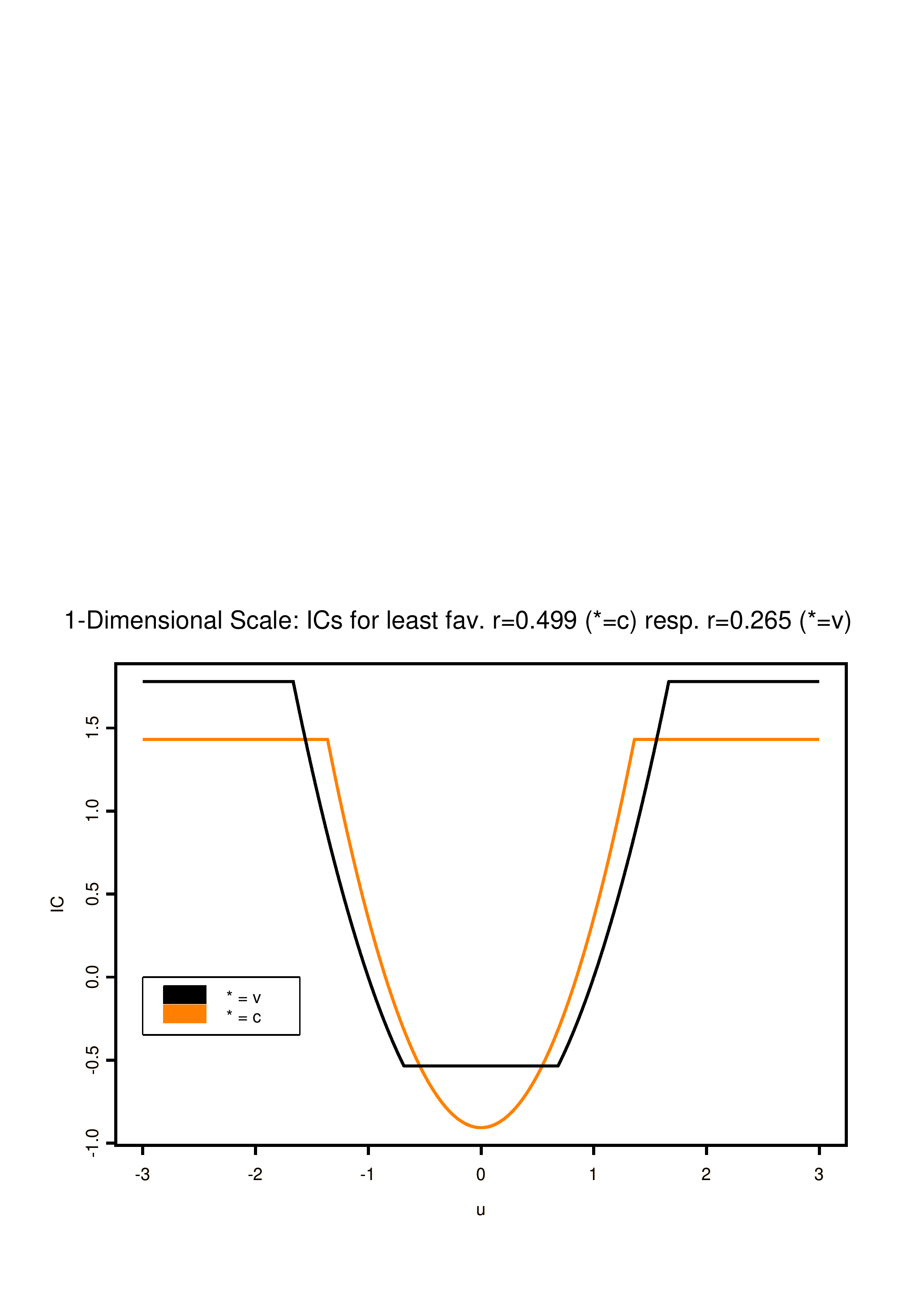}
\end{center}
\end{frame}
\begin{frame}\vspace{-5cm}\begin{center}
\includegraphics[width=0.8\linewidth]{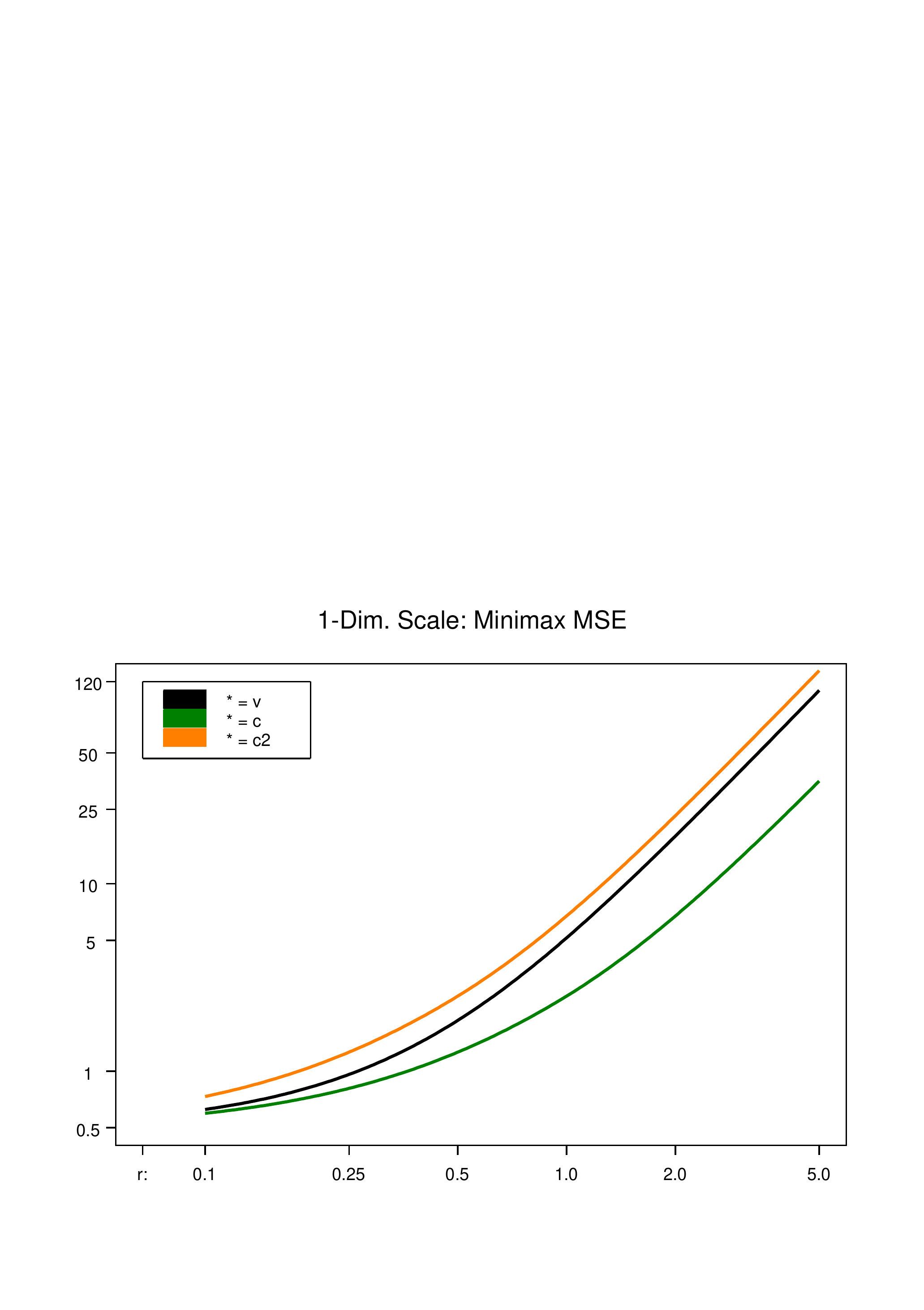}
\end{center}
\end{frame}
\begin{frame}\vspace{-5cm}\begin{center}
\includegraphics[width=0.8\linewidth]{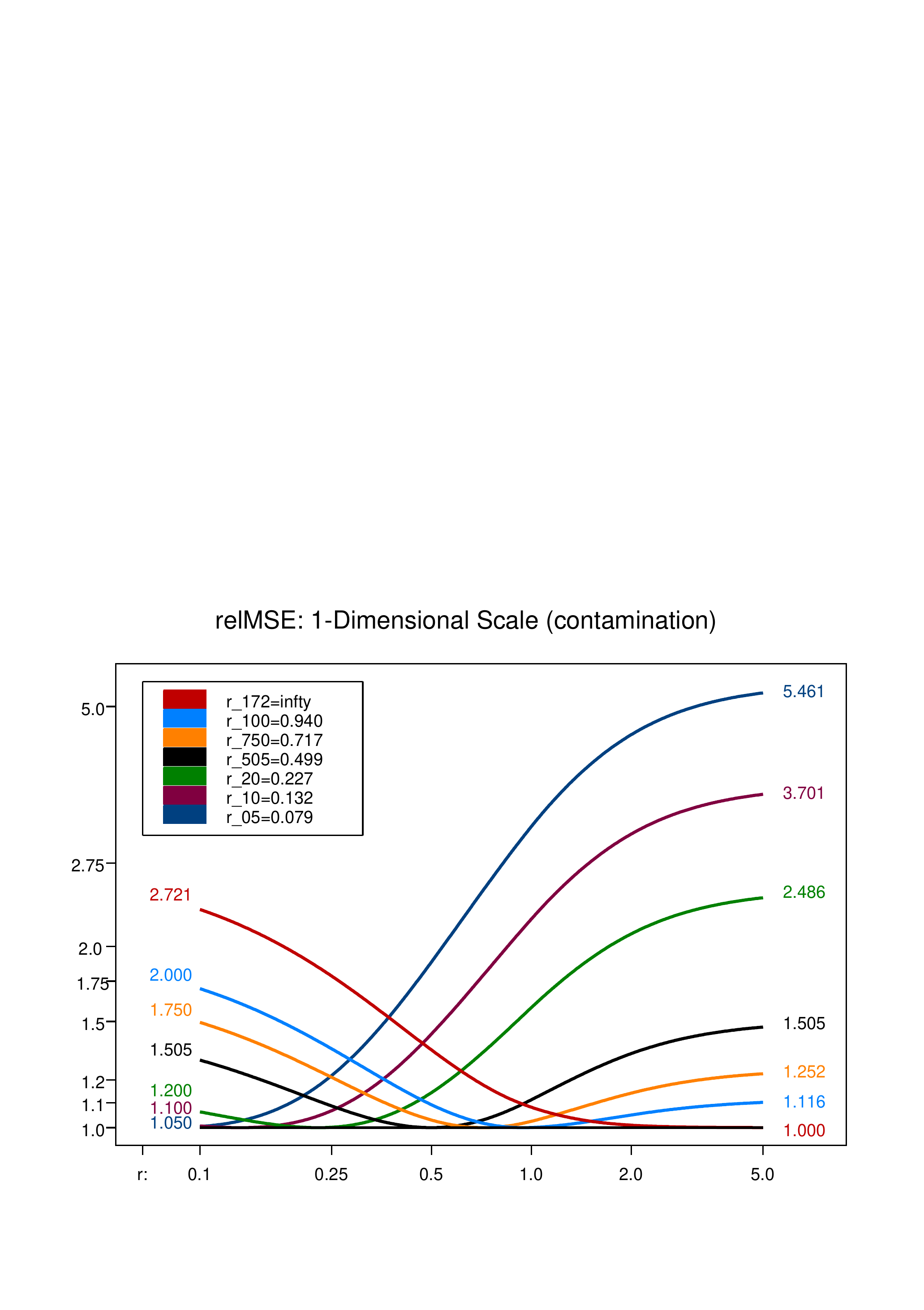}
\end{center}
\end{frame}
\begin{frame}\vspace{-5cm}\begin{center}
\includegraphics[width=0.8\linewidth]{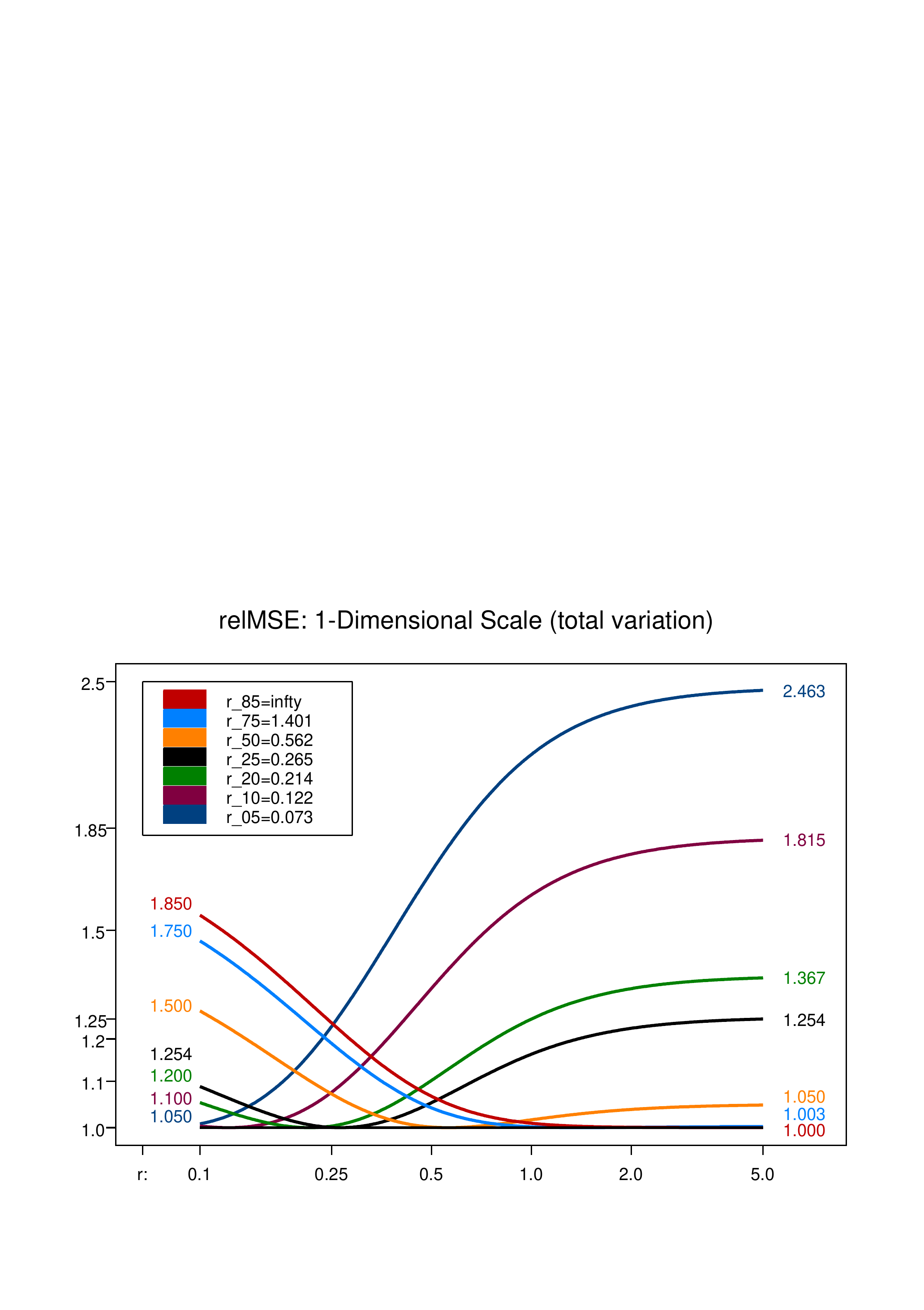}
\end{center}
\end{frame}
\subsection{6.8 Summary}\begin{frame}
\blu{\large\bf Summary}
\par \bigskip
1) Estimation of the unknown radius hardly pays, provided one employs the
radius-minimax estimator. The increase of its risk with respect to the
radius-optimal procedure is moderate to small.
\par
In all our models, it is $\le 12.5$\%, if the radius may be specified
to belong to some interval $ [\hskip6\ei \frac{1}{3}r,3r]$ for any~$r$.
\par \medskip
2) The minimax radii are small: 
                                5--6\% contamination, at sample size 100. 
\par \medskip
3) The radius-minimax estimator for completely unknown radius stays the
same for a variety of convex risks which are homogeneous in bias
and (square root) variance; e.g., $L_p$-loss, confidence levels.
\hfill \mbox{\tiny Ruckdeschel, Rieder (2004)}
\par \hfill \hfill \mbox{\tiny Rieder, Kohl Ruckdeschel (2008)}
\par \null \vfill
\end{frame}
\section{References}
\begin{frame}\null \frenchspacing \tiny
Beran, R.J. (1974):\hskip.5em Asymptotically efficient and
       adaptive rank estimates in location models.
       {\sl Ann. Statist.\/} {\bf 2} 63--74.
\par \vskip.125pt plus .125pt
\mbox{Beran,~R.J.~(1976)}:\enskip
         Adaptive estimates for autoregressive processes.
        {\sl Ann. Inst. Statist. Math.\/}~{\bf 11} \mbox{432--452}.
\par \vskip.125pt plus .125pt
 \mbox{Bickel, P.J.~(1982)}:\hskip.5em On adaptive estimation.
   {\sl Ann. Statist.\/} {\bf 10} 647--671.
\par \vskip.125pt plus .125pt 
 \mbox{Bickel, P.J.} et al.~(1993):\hskip.5em
  {\sl Efficient and Adaptive Estimation for Semiparametric Models\/}.
  Springer, New York.
\par \vskip.125pt plus .125pt
 Birg\'e, L. (1980): Approximation dans les espaces m\'etrique
   et th\'eorie de l'estimation. Ph.D.~Thesis, University of Paris.
\par \vskip.125pt plus .125pt
          Drost, F.C., Klaassen, C.A.J., Werker, B.J.M. (1997):
          Adaptive estimation in time-series models.
         {\sl Ann. Statist.\/} {\bf 25} 786--817.
\par \vskip.125pt plus .125pt
          Fischer, J. (2006): {\sl Robuste Sch\"atzung im semiparametrischen
          Mixture-Modell\/}. Diploma thesis, U Bayreuth.
\par \vskip.125pt plus .125pt
         Hampel, F.R. et al.~(1986):\hskip.5em
         {\sl \mbox{Robust} Statistics---The Approach Based on
         Influence Functions\/}. Wiley, New York.
\par \vskip.125pt plus .125pt
 \mbox{Huber,~P.J.~(1981)}:\hskip.5em
           {\sl Robust Statistics\/}. Wiley, New York.
\par \vskip.125pt plus .125pt
 \mbox{Huber,~P.J.} and Strassen,~V.(1973):\hskip.5em
        Minimax tests and the Neyman--Pearson lemma for capacities.
        {\sl Ann.\ Statist.\/}~{\bf1} \mbox{251--263}.
\par \vskip.125pt plus .125pt
          Huber--Carol, C.~(1970):\hskip.5em
    {\sl \raisebox{0ex}[1ex][0ex]{\'Etude} asymptotique de tests robustes\/}.
         Th\`ese de \mbox{Doctorat}, ETH~Z\"urich.
\par \vskip.125pt plus .125pt
    Kohl, M. (2005):
    {\em Numerical Contributions to the Asymptotic Theory of Robustness}.
    PhD thesis. University of Bayreuth.
\par \vskip.125pt plus .125pt
        Kohl, M. and Ruckdeschel, P. (2008):
       ROptEst: Opt. robust estimation.
      {\em R ver. 0.6.3}. 
      URL http://robast.r-forge.r-project.org
\par \vskip.125pt plus .125pt
       Kreiss, J.-P.~(1987): On adaptive estimation in stationary
       ARMA processes. {\sl Ann.\ Statist.\/}~{\bf15} \mbox{112--133}.
\par \vskip.125pt plus .125pt  
       Pfanzagl, J., Wefelmeyer, W. (1982): {\sl
       Contributions to a General Asymptotic Statistical Theory\/}.
       Springer LN in Statistics \#13. 
\par \vskip.125pt plus .125pt
  \mbox{Rieder,~H.~(1977)}:\hskip.5em
         Least favorable pairs for special capacities.
         {\sl Ann.\ Statist.\/} {\bf5} \mbox{909--921}.
\par \vskip.125pt plus .125pt
        \mbox{Rieder,~H.~(1978)}:\hskip.5em
         A robust asymptotic testing model.
         {\sl Ann.\ Statist.\/}~{\bf6} \mbox{1080--1094}.
\par \vskip.125pt plus .125pt                 
 \mbox{Rieder,~H.~(1994)}:\hskip.5em
           {\sl Robust Asymptotic Statistics\/}. Springer, New York.
\par \vskip.125pt plus .125pt
   \mbox{Rieder,~H.~(2000)}:\hskip.5em
  Neighborhoods as nuisance parameters? Robustness vs.\ semiparametrics.
         Discussion Paper Nr. 25, SFB 373. 
\par \vskip.125pt plus .125pt \mbox{Rieder,~H.~(2000)}:\hskip.5em
       One-sided confidence about functionals over tangent cones.
         Discussion Paper Nr. 26, SFB 373. 
\par \vskip.125pt plus .125pt
        Rieder, H., Kohl, M. and Ruckdeschel, P. (2008):\hskip.5em
        The cost of not knowing the radius.
        {\sl Stat. Meth.\& Appl.} {\bf 17} 13--40.
\par \vskip.125pt plus .125pt
         Ruckdeschel, P. and Rieder, H. (2004):
         Optimal influence curves for general loss functions.
         {\sl Statistics\&Decisions\/}~{\bf 22} 201–223.
\par      Ruckdeschel, P. (2006):
          A motivation for $1/\!\sqrt{n\,}$-shrinking neighborhoods.
          {\sl Metrika\/} {\bf 63} 295–-307.
\par \vskip.125pt plus .125pt
  Ruckdeschel, P., Hable, R., Rieder, H. (2010):
          Optimal robust ICs in semiparametric regression.
           {\sl JSPI\/} {\bf 140} 226--245.
\par \vskip.125pt plus .125pt
          Ruckdeschel, Kohl, Rieder (2010):
          Infinitesimally Robust estimation in general smooth models. 
          {\sl Stat. Meth.\& Appl.} {\bf 19} 333--354.
\par \vskip.125pt plus .125pt
         Ruckdeschel, P. and Rieder, H. (2010):
         Fisher information of scale.
         {\sl Statistics and Probability Letters\/}~{\bf 80} 1881-–1885.
\par \vskip.125pt plus .125pt
         \mbox{Shen, L.Z.~(1994)}:\hskip.5em
   Optimal robust estimates for semiparametric symmetric location models.
   {\sl Statistics\&Decisions\/}~{\bf 12} 113--124.
\par \vskip.125pt plus .125pt
      \mbox{Shen, L.Z.~(1995)}:\hskip.5em
   On optimal B-robust influence functions in semiparametric models.
   {\sl Ann. Statist.\/}~{\bf 23} 968--989.
\par \vskip.125pt plus .125pt
   Stabla, Th. (2005): {\sl Robuste adaptive Sch\"atzung\/}.
   Diploma thesis, U Bayreuth.
\par \vskip.125pt plus .125pt
   Stein, C. (1956):\hskip.5em
  {\em Efficient nonparametric estimation and testing}.
  In: Proceedings of the Third Berkeley Symposium on Mathematical Statistics
  and Probability, 1954--1955, vol. I , 187--195.
  UC Press, Berkeley and Los Angeles.
\par \vskip.125pt plus .125pt
 Stone, C.~(1975):\hskip.5em
   Adaptive maximum likelihood estimation for a location parameter.
   {\sl Ann. Statist.\/} {\bf 3} 267--284.
\par \vskip.125pt plus .125pt
 van der Vaart, A.W.~(1998):\hskip.5em
   {\sl Asymptotic Statistics\/}. CUP, Cambridge.
\par \null\vfill \nonfrenchspacing \nz \end{frame}
\end{document}